\documentclass[a4paper,leqno,11pt]{amsart}

\usepackage[english]{babel}

\makeatletter
\@namedef{subjclassname@2020}{\textup{2020} Mathematics Subject Classification}
\makeatother

\usepackage{amsfonts}
\usepackage{amsmath}
\usepackage{amsthm}
\usepackage{amssymb}
\usepackage{indentfirst}
\usepackage{color}
\usepackage{tabularx}
\usepackage{graphicx}
\usepackage{esint}
\allowdisplaybreaks

\usepackage[a4paper,top=3.5cm,bottom=1cm,left=2.5cm,right=2.5cm]{geometry}

\usepackage{subcaption}

\usepackage{mathrsfs}
\usepackage{mathtools}
\usepackage{enumitem}
\usepackage[normalem]{ulem}

\definecolor{blue_links}{RGB}{13,0,180} 

\usepackage{hyperref}
\hypersetup{
    colorlinks=true, 
    linktoc=all,     
    linkcolor=blue_links,  
    citecolor=blue_links,
    urlcolor=blue_links,
}

\usepackage[utf8]{inputenc}
\usepackage[T1]{fontenc}
\usepackage[english]{babel}

\numberwithin{equation}{section}

\theoremstyle{plain}
\begingroup
\theoremstyle{plain}
\newtheorem{theorem}{Theorem}[section]
\newtheorem{corollary}[theorem]{Corollary}
\newtheorem{proposition}[theorem]{Proposition}
\newtheorem{lemma}[theorem]{Lemma}
\newtheorem{definition}[theorem]{Definition}

\theoremstyle{remark}
\newtheorem{remark}[theorem]{Remark}
\endgroup

\theoremstyle{definition}
\theoremstyle{remark}

\mathchardef\emptyset="001F

\newcommand{\strain}{\boldsymbol{\varepsilon}} 

\newcommand{\R}{\mathbb{R}}
\newcommand{\E}{\mathcal{E}}

\newcommand{\dive}{\mathrm{div}}

\newcommand{\C}{\mathbb{C}}
\newcommand{\e}{\mathrm{E}}

\newcommand{\G}{\mathcal{G}}

\newcommand{\di}{\mathrm{d}}

\newcommand{\A}{\mathcal{A}}
\newcommand{\HH}{\mathcal{H}}
\newcommand{\Om}{\Omega}

\newcommand{\M}{\mathbb{M}}

\newcommand{\F}{\mathcal{F}}

\newcommand{\J}{\mathcal{J}}

\newcommand{\argmin}{\mathrm{argmin}}

\definecolor{dred}{rgb}{.8,0,0}
\definecolor{ddmagenta}{rgb}{0.7,0,0.9}
\definecolor{ddcyan}{rgb}{0,0.2,1.0}
\definecolor{Orchid}{rgb}{0.7,0.4,0}

\newcommand{\eps}{\varepsilon}

\newcommand{\V}{\mathrm{V}}

\newcommand{\coloneq }{\hspace{1pt}\raisebox{0.74pt}{\scalebox{0.8}{:}}\hspace{-2.2pt}=}

\title[Topology optimization with microstructure]{Phase-field topology
  optimization \\ with periodic microstructure}
\author[S. Almi]{Stefano Almi}
\address[Stefano Almi]{Department of Mathematics and Applications {\it
    R.~Caccioppoli}, via Cintia, Monte S.~Angelo, 80126 Naples, Italy.}
\email{stefano.almi@unina.it}
\urladdr{https://www.docenti.unina.it/stefano.almi}

\author[U. Stefanelli]{Ulisse Stefanelli} 
\address[Ulisse Stefanelli]{Faculty of Mathematics, University of
  Vienna, Oskar-Morgenstern-Platz 1, A-1090 Vienna, Austria,
Vienna Research Platform on Accelerating
  Photoreaction Discovery, University of Vienna, W\"ahringerstra\ss e 17, 1090 Wien, Austria,
 \& Istituto di
  Matematica Applicata e Tecnologie Informatiche {\it E.~Magenes}, via
  Ferrata 1, I-27100 Pavia, Italy
}
\email{ulisse.stefanelli@univie.ac.at}
\urladdr{http://www.mat.univie.ac.at/$\sim$stefanelli}

\date{\today}
 \subjclass[2020]{74Q05, 
 			   74P10,  
			   49Q10,  	
			   49J20, 	
			   49K20,  	
			   }
\keywords{Periodic microstructure, topology optimization,
  homogenization, sharp-interface limit, first-order
conditions.}

\begin{document}
  \begin{abstract}

Progresses in additive manufacturing technologies allow the realization
of finely graded microstructured materials with tunable mechanical
properties. This paves the way to a wealth of innovative applications,
calling for the combined design of the macroscopic mechanical piece
and its underlying microstructure. In this context, we investigate a topology optimization problem for an elastic medium
featuring a periodic microstructure.  The optimization problem is variationally
formulated as a bilevel minimization of phase-field type. By resorting
to $\Gamma$-convergence techniques, we characterize the homogenized
problem and investigate the corresponding sharp-interface limit. 
First-order optimality conditions are derived, both at the
homogenized phase-field
and at the sharp-interface level. 
\end{abstract}
\maketitle

\section{Introduction}
\label{s:intro}

Additive manufacturing  allows the direct  production of 
three-dimensional (3D) components  by successive layer
deposition.  The macroscopic shape of the manufactured
 workpiece has to be optimized with respect to its mechanical
 function. This is the general task of {\it Topology Optimization}
 (TO), which 
applies to a number of different shape design problems, from
mechanical engineering, to aerospace and 
automotive, to architectural engineering, to biomechanics
\cite{MR2008524}.

Recent technological advances allow the realization of microscopic
patterns of graded materials within a single macroscopic piece \cite{Carraturo}. As a
consequence, the material response can be modulated by designing this
microstructure, which calls for an optimization process at the
microscale, as well.
From the computational viewpoint, combined microstructure and topology
optimization problems have been recently investigated in a
number of different settings. The reader is referred to
\cite{Barbarosie,Christensen,Geoffroy-Donders,Harbrecht,Kim,Valentin} among
others, as well as to the
recent survey \cite{Wu}. 

The microstructure and topology-optimization problem is recasted in
\cite{Allaire19, Allaire20} within the frame of the {\it homogenization method}
\cite{MR1859696}. In 2D and 3D, admissible microstructures are periodic perforations
by cubic holes,  possibly  reoriented via a given rotation
field. The macroscopic effect of the microstructure is resolved
by computing the homogeneized elastic response. Then, the optimal macroscopic
shape is numerically investigated.

Manufacturing constraints are additionally considered in \cite{Conti}. Given an admissible macroscopic response, a
corresponding microstructure is identified by
optimizing a microscopic criterion, which penalizes microstructural volume and surface in
the periodicity cell. A regularization is then obtained by a
phase-field approach, allowing for numerical treatment, both at the
micro- and at the macroscale.

In this paper, we present the theoretical analysis of a
 topology-optimization problem in
presence of a periodic microstructure.  We focus on a simplified, reference
setting, amenable to a rigorous mathematical treatment.  Similar to~\cite{Conti} (see, e.g.,~\cite{Blank2, Garcke-23, Bourdin, Burger, Penzler} for further applications in elasticity), our analysis builds upon a phase-field approach.  Shape and microstructure are described via two scalar
fields. At first, the actual distribution of
material within an a-priori given design region $\Omega  \subseteq  
{\mathbb R}^n$ is described by a level set of the phase field
 $\varphi\colon \Omega \to [0,1]$.  The actual optimal topology to be
identified corresponds to the set $\{\varphi=1\}$. As it
customary in topology optimization, one considers the whole
design domain $\Om$ to be filled by an elastic medium, where the
region $\{\varphi<1\}$ has to be interpreted as a very compliant {\it
  Ersatz} material. The integral of $\varphi$ on $\Om$ is
constrained in order to exclude the trivial case $\varphi\equiv 1$. 

The microstructure of the medium is described via the periodic
function  $m \colon Y \to [1,2]$  defined on the unit cube $Y:=[0,1)^n$. The
value of $m$ represents the grading between two
possible constituents of the microstructure. The mean of the microstructure variable is 
constrained in $(1,2)$, so to avoid trivial solutions. By 
assuming $\eps>0$ to be the microstructural scale, the 
elastic response of the medium is then encoded in the choice of the
scale-dependent elastic tensor
$$\C_\eps(\varphi,m)(x) \coloneq \varphi(x)m\left(\frac{x}{\eps} \right)\C_1 +
(1{-}\varphi(x))\C_2.$$
Here, the elastic tensor $\C_1$ describes the elastic
response of the medium (namely, for $\varphi=1$) for constituent $m=1$, whereas the different
behavior in case of constituent $m=2$ is simply modeled by $2\C_1$. On the other
hand, $\C_2$ is the elastic tensor of the very soft Ersatz material. Both $\C_1$ and $\C_2$ are
asked to be major and minor-symmetric and positive definite, possibly
with $\C_2<\!<\C_1$. 

By postponing all the necessary detail to the forthcoming Sections, for each suitably given pair $(\varphi,m)$, the elastic equilibrium of
the body is described via the displacement  $u\colon \Om \to \R^n$,  solving
the equilibrium system
$$-\dive \Big(\C_\eps(\varphi,m)\strain(u) \Big) = \varphi f \ \
\text{in} \ \Om, \quad \C_\eps(\varphi,m)\strain(u) \nu = g \ \
\text{on} \ \Gamma_N, \quad u=0 \ \ \text{on} \ \Gamma_D$$
where $\strain(u) $ stands for the symmetrized strain $\strain(u) =
(\nabla u +\nabla^\top u)/2$, $f$ encodes body forces, $g$ is the boundary
traction at $\Gamma_N  \subseteq  \partial \Om$, $\nu$ is outward normal
to the boundary $ \partial \Om$,
and the body is clamped ($u=0$) at $\Gamma_D \subset \Om$. 
The latter
equilibrium problem, can be formulated variationally, by identifying
the displacement
$u$ as the unique minimizer of an elastic energy
$\E_\eps(\varphi,m,\cdot)$, see \eqref{e:energy} below. Given the
triple $(\varphi,m,u)$, our optimization aim is to minimize the  functional
$$\J_{\eps} (\varphi, m, u) \coloneq  \mathcal{C}(\varphi,u)  + \frac12 \fint_{\Om} \Big| \nabla m
\Big(\frac{x}{\eps} \Big) \Big|^{2} \, \di x + \frac12 \int_{\Om} |
\nabla \varphi|^{2} \, \di x \,$$
where $\mathcal{C}(\varphi,u)$ is the classical {\it compliance}
$$\mathcal{C}(\varphi,m)  \coloneq \int_{\Om}
\varphi f{\, \cdot\,} u \, \di x + \int_{\Gamma_{N}} g{\, \cdot\,} u
\, \di \HH^{n-1}$$
($\HH^{n-1}$ denoted the $(n-1)$-dimensional Hausdorff (surface) measure) and the
gradient terms instead penalize microstructural interfaces at the
microlevel and material interfaces at the macrolevel,
respectively. Our departing point will be the study of the

$$\boxed{\begin{array}{l}
  \text{Microstructure-Topology Optimization (MTO) problem:}\\[2mm]
{}\qquad\min_{(\varphi,m)} \Big\{ \J_\eps(\varphi,m,u) \: : \: u =\argmin\,
  \E_\eps(\varphi,m,\cdot) \Big\}\
       \end{array}
     }
     $$
     \smallskip
     
\noindent The latter is proved to admit a solution in Proposition
\ref{prop:ex:MTO}.

A first focus of our investigation is on the $\eps\to 0$
{\it homogenization} limit of (MTO). As the scale of the microstructure goes to $0$,
the effective material response can be described by an homogenized
elastic tensor $\C^*(\varphi,m)$, which is completely characterized in
Section \ref{s:hom} in terms of auxiliary cell problems. This in
particular corresponds to a homogenized equilibrium problem, where
$\C_\eps$ is replaced by $\C^*$, as well as a 
homogenized elastic energy $\E(\varphi,m,u)$, see \eqref{e:energy-2}
below. We investigate the homogenized cost functional $\J$ taking the form 
$$\J (\varphi, m, u) \coloneq \mathcal{C}(\varphi,u)+ \frac12
\int_{Y}| \nabla m |^{2} \, \di y + \frac12 \int_{\Om} | \nabla
\varphi|^{2} \, \di x \,,$$
where, compared with $\J_\eps$, variations of the microstructure
variable $m$ are simply evaluated in the microscopic variable. Correspondingly, we consider the 

  $$\boxed{\begin{array}{l}
  \text{Homogenized Microstructure-Topology Optimization (HMTO) problem:}\\[2mm]
\qquad\qquad\quad \min_{(\varphi,m)} \Big\{ \J(\varphi,m,u) \: : \: u =\argmin\,
  \E(\varphi,m,\cdot) \Big\} 
    \end{array}}
    $$
       \smallskip
     
\noindent 
This is proved  to  admit solutions in Proposition \ref{p:existence}. In
particular, we link the (HMTO) problem to the (MTO) problem for $\eps>0$
via $\Gamma$-convergence in Proposition \ref{p:Gamma-J}. In addition,
we are able to provide first-order optimality conditions for solutions to
(HMTO) in  Theorem
\ref{t:optimality-deterministic}. Here, admissible variations are
required to respect the volume constraints imposed to the phase field
$\varphi$ and the microstructure $m$, see Definition \ref{d:a} below. 

We then turn to the consideration of the sharp-interface limit of the
(HMTO) problem. This ensues from considering a Modica-Mortola generalization
of the cost functional $\J_\eps$ as
\begin{align*} 
\mathcal{J}_{\varepsilon}^{ \text{s}} (\varphi, m, u) &\coloneq 
                                                     \mathcal{C}(\varphi,u)+
                                                     \frac{\eps}{2}\int_{\Om}
                                                     |
                                                     \nabla
                                                     \varphi|^{2} \di
                                                     x + \int_{\Om}
                                                     \frac{1}{2\eps}H(\varphi)
\, \di x \\
&+ \frac{\eps}{2}\fint_{\Om} | \nabla m| ^{2} \Big(
                                                     \frac{x}{\eps}
                                                     \Big) \di x+
                                                     \fint_{\Om} \frac{1}{2\eps} H \Big(m\Big(
                                                     \frac{x}{\eps}
                                                     \Big) -1
                                                     \Big) \di x \nonumber
\end{align*} 
where the function  $H \colon \R \to [0,+ \infty)$ is nonnegative and $H(r)=0$
iff $r = 0$ or $r=1$.  By taking the limit $\eps \to 0$, one identifies
the sharp-interface functional
\begin{align*}
\mathcal{J}^{ \text{s}} (\varphi, m, u) \coloneq  \mathcal{C}(\varphi,u) +
  c_{H} \mathcal{H}^{n-1}(\{\varphi = 1\}\cap \Om) +  c_{H}
  \mathcal{H}^{n-1} (\{ m=2\}\cap Y) \big)\,,
\end{align*}
featuring the lengths of the interfaces of the regions $\{\varphi =
1\}$ and $\{m = 2\}$ in $\Omega$ and $Y$, respectively. Here, the
constant $c_H>0$ is directly computed from $H$. We hence
consider the

$$
\boxed{\begin{array}{ll}
  \text{Sharp-interface Microstructure-Topology Optimization (SMTO) problem:}\\[2mm]
\qquad\quad \min_{(\varphi,m)} \Big\{ J^{ \text{s}}(\varphi,m,u) \: : \: u =\argmin\,
  \E(\varphi,m,\cdot) \Big\}
       \end{array}}
     $$\smallskip

     \noindent
This is proved to admit a solution in Proposition \ref{p:sharp-conv}
by means of a $\Gamma$-convergence argument. First-order optimality
conditions for the (SMTO) problem are then provided by Theorem
\ref{t:optimality}. Compared with the optimality conditions for the
(HMTO) problem, the situation is here more delicate, for admissible variations
have to be additionally adapted to the identification of
interfaces, see Definitions \ref{d:variations}-\ref{d:variations2}
below.

Before closing this introduction, we would like to mention a suite of
aspects, which are currently not covered by our analysis, but that we
deem as relevant for future work.

At first, let us remark that our choice for the possible
microstructure of the material would call for extensions in
different directions. The exact periodicity assumption, although
convenient from the mathematical viewpoint and well justified in some
contexts, may be less suited for some applications. It would be
relevant to consider quasiperiodic or even stochastic settings
instead.  The actual technical manufacturability could also be
considered to limit the choice of  the possible microstructures.  

In addition, it would be advantageous to tackle the possibility of
changing the orientation of the microstructure in various regions of
the elastic piece. In 3D printing such possibility would be realizable
by specifying deposition patterns, see \cite{Allaire19,Allaire20}.
More generally, it would be convenient to possibly change the
microstructure in different regions of the microscale.

Eventually, the underlying mechanical setting might be adapted to
different situations, including alternative mechanical responses
(inelastic, dynamic, thermalized, among others)  \cite{A-S-21, A-S-22, Boissier, Desai, Maury},  different
optimization goals (stress and strain control, stiffness maximization,
natural frequencies and damping design, vibration stability, etc.)  \cite{Auricchio, Joubert} 
under other constraints (material availability, weight, or
cost, manufacturing constraints, structural requirements, for
instance)  \cite{Conti, Garcke-23-2}. 

The paper is organized as follows. We present the mechanical and the
functional setting in detail in Section \ref{s:det} where the
existence for the (MTO) problem is also discussed. Section \ref{s:hom}
is then devoted to the (HMTO) problem. In particular, the first-order
optimality conditions are presented in Subsection
\ref{s:optimality}. Eventually, the (SMTO) problem is discussed in
Section \ref{s:sharp} and the corresponding first-order optimality
conditions are in Subsection \ref{s:opt-sharp}.

\section{The Microstructure-Topology Optimization problem}
\label{s:det}

 We devote this section to the specification of the mechanical
model and of its variational setting. 

Let $\Om\subseteq \R^{n}$ be open, bounded, with Lipschitz
boundary~$\partial \Om$, let $\Gamma_{D}, \Gamma_{N}
\subseteq \partial \Om$  be open in the topology of $\Omega$ and
  be such that $\HH^{n-1}(\Gamma_{D}) >0$ and
$\overline{\Gamma}_{D} \cap \overline{\Gamma}_{N} = \emptyset$,  
where we recall that $\HH^{n-1}$ is the $(n-1)$-dimensional Hausdorff measure. We 
fix the applied  density per unit $\varphi$ of  volume forces
$f \in L^{2}(\Om; \R^{n})$ and boundary-traction  density   $g \in L^{2}(\Gamma_{N}; \R^{n})$. We define
\begin{align*}
 H^{1}_{\Gamma_{D}} (\Om; \R^{n}) \coloneq \{ u \in H^{1}(\Om; \R^{n}) : \, u= 0 \text{ on $\Gamma_{D}$}\}\,,
\end{align*}
denote by $Y$ the unit  cube  $[0, 1)^{n}$, and  let  $H^{1}_{\#} (Y;
 A)$  be  the set of $H^{1}$- and $Y$-periodic functions with
values in  the set~$A  \subseteq  \R^{m}$.  We further fix $\C_{1}$
and $\C_{2}$ the elasticity tensors of the material to be distributed
and  of the Ersatz material.   Both are supposed to be 
major and minor symmetric and  positive definite, i.e., there exist $0<\alpha \leq \beta < \infty$ such that
\begin{equation}\label{e:hp-C}
\alpha | \e |^{2} \leq \C_{i} \e {\, \cdot\,} \e \leq \beta | \e|^{2}
\qquad \text{for  all  $\e \in \M^{n}_{S}$ and $i = 1, 2$\,,}
\end{equation}
where $\M^{n}_{S}$ denotes the space of symmetric matrices of order
$n$.  We recall that the {\it compliance}  $\mathcal{C} \colon 
H^{1}(\Om) \times   H^{1}(\Om; \R^{n})$  is
defined as
$$\mathcal{C}(\varphi,u)\coloneq \int_{\Om}
\varphi f{\, \cdot\,} u \, \di x + \int_{\Gamma_{N}} g{\, \cdot\,} u
\, \di \HH^{n-1}.$$

For every $\varepsilon>0$, every $u \in H^{1}(\Om; \R^{n})$, 
every $\varphi \in H^{1}(\Om)$, and every $m \in H^{1}_{\#}(Y)$  we define the energy functional
\begin{align}
\label{e:energy}
\E_{\varepsilon} ( \varphi,m, u) \coloneq &\ \frac12\int_{\Om} \left( \varphi (x) m\left(\frac{x}{\eps}\right) \C_{1} + (1 {-} \varphi(x) ) \C_{2}\right) \strain(u) {\, \cdot\,} \strain(u) \, \di x 
-   \mathcal{C}(\varphi,u) 
\end{align}
where  we recall that  $\strain(u)  \coloneq ({\nabla u + \nabla^{\top} u})/{2}$ denotes the symmetric part of the gradient of~$u$. To
 streamline  the notation, we set
\begin{displaymath}
\C(\varphi, m) ( x, y ) \coloneq  \varphi (x) m(y) \C_{1} + (1 {-} \varphi(x) ) \C_{2} \qquad \text{for $(x, y) \in \Om \times Y$.}
\end{displaymath}
With  a slight  abuse of notation, we still denote by $\C(\varphi, m) ( x,
\frac{x}{\varepsilon})$ the tensor used in~\eqref{e:energy}.

 In
\eqref{e:energy} the  phase-field  function $\varphi$ 
indicates the relative concentration of the
microstructured material with respect to Ersatz material.  The function $m (x/\varepsilon)$  defines the
  the microstructure at scale~$\varepsilon$,  which is assumed
to affect the elastic response by (multiplicatively) influencing the
elasticity tensor~$\C_{1}$.  In what follows, $\varphi$ takes values in
the interval~$[0, 1]$, while~$m$ takes values in the~$[1, 2]$, so that
 for all~$(\varphi,m) \in H^1(\Omega;[0,1]) \times H^1_\#(Y;[1,2]) $
 the energy~$u \in H^1_{\Gamma_D}(\Omega;\R^n)\mapsto {\mathcal
  E}_\varepsilon(\varphi,m,u)$ from~\eqref{e:energy}  is coercive.

As in many applications one aims at designing light structures or at
reducing productions costs,  in the following we impose some
constraints on phase field~$\varphi $ and on the microstructure~$m$,
having also the effect of excluding
trivial solutions.   Namely, we fix $V \in  (0,
\mathcal{L}^{n}(\Om))$  and $W \in  (1,2)$  and define,
for $\varepsilon>0$, the set of admissible  phase fields and
microstructures  as
\begin{displaymath}
 \A_{\varepsilon} \coloneq \left\{ (\varphi , m) \in H^{1}(\Om; [0, 1]) \times H^{1}_{\#} (Y; [1, 2]): \, \int_{\Om} \varphi \, \di x \leq V \text{ and } \fint_{\Om} m\Big(\frac{x}{\varepsilon} \Big) \, \di x \leq W\right\}\,.
\end{displaymath}

For every $u \in H^{1}_{\Gamma_{D}}(\Om; \R^{n})$ and every $(\varphi,
m) \in  \A_{\eps}$ we define the $\eps$-cost functional
\begin{align}
\label{e:eps-cost}
\J_{\eps} (\varphi, m, u) \coloneq   & \    \mathcal{C}(\varphi,u)
+ \frac12 \fint_{\Om} \Big| \nabla m \Big(\frac{x}{\eps} \Big) \Big|^{2} \, \di x + \frac12 \int_{\Om} | \nabla \varphi|^{2} \, \di x \,,
\end{align}
and set $\J_{\eps} (\varphi, m, u) =  \infty$ otherwise in $H^{1}_{\Gamma_{D}}(\Om; [0,1]) \times H^{1}_{\#}(Y; [1,2]) \times H^{1}(\Om; \R^{n})$. 

%
 
For $\eps>0$, we consider the {\it microstructure-topology
  optimization}  (MTO)  problem

\begin{align}
\label{e:pb-eps}
& \min_{\substack{\varphi \in H^{1}(\Om; [0, 1]) ,\\ m \in
  H^{1}_{\#}(Y; [ 1, 2])} } \, \J_{\eps}(\varphi, m, u)  \ \  \ \
  \text{subject to } \min_{u \in H^{1}_{\Gamma_{D}} (\Om; \R^{n}) } \,
  \E_{\eps}(\varphi, m, u) \,. \tag{MTO}
\end{align}  

For later use, we denote by~$S_{\eps} \colon L^{1}(\Om; [0, 1]) \times  L^{1}_{\#}(Y; [ 1, 2]) \to H^{1}_{\Gamma_{D}} (\Om; \R^{n})$ the map that to~$(\varphi, m) \in L^{1}(\Om; [0, 1]) \times L^{1}_{\#}(Y; [ 1, 2])$ associates the unique solution to the elastic equilibrium problem in~\eqref{e:pb-eps}.

In the following proposition we show existence of solutions to~\eqref{e:pb-eps}.

\begin{proposition}[Existence for the \eqref{e:pb-eps} problem]
\label{prop:ex:MTO}
 For all $\varepsilon > 0$ problem~\eqref{e:pb-eps} admits a
solution. 
\end{proposition}

\begin{proof}
Let $(\varphi_{j}, m_{j}, u_{j}) \in H^{1}(\Om; [0,1]) \times
H^{1}_{\#}(Y; [1, 2]) \times H^{1}_{\Gamma_{D}} (\Om; \R^{n})$ be a
minimizing sequence for~\eqref{e:pb-eps}. Then,  
$(\varphi_{j}, m_{j}) \in \A_{\eps}$  for every $j$ and that, up
to a subsequence, $\varphi_{j} \rightharpoonup \varphi$ weakly in
$H^{1}(\Om)$ and $u_{j} \rightharpoonup u$ weakly in
$H^{1}_{\Gamma_{D}} (\Om; \R^{n})$. Moreover, $\varphi \in H^{1}(\Om;
[0,1])$. Let us show that $m_{j}$ is bounded in $H^{1}_{\#}(Y)$. For
every $\eps>0$ and every $z \in \mathbb{Z}^{n}$, we set $Q^{\eps}_{z}
\coloneq \varepsilon z + (0, \varepsilon)^{n}$. For every open subset
$ \emptyset \not =  \Om' \Subset \Om$ we define $\mathbb{Z}^{n}_{\varepsilon, \Om'} \coloneq \{ z \in \mathbb{Z}^{n}: \, Q^{\eps}_{z} \cap \Om' \neq \emptyset\}$. Then, up to taking a smaller~$ \emptyset \not = \Om' \Subset \Om$ we may assume that $Q^{\eps}_{z} \subseteq \Om$ for every $z \in \mathbb{Z}^{n}_{\eps, \Om'}$. Hence, we have that
\begin{align}\label{e:12}
\fint_{\Om} \Big| \nabla m_{j}\Big(\frac{x}{\eps} \Big) \Big|^{2}\, \di x & \geq \frac{1}{\mathcal{L}^{n}(\Om)} \sum_{z \in \mathbb{Z}^{n}_{\varepsilon, \Om'}} \int_{Q^{\eps}_{z}} \Big| \nabla m_{j} \Big(\frac{x}{\eps} \Big) \Big|^{2}\, \di x 
\\
&
=  \frac{1}{\mathcal{L}^{n}(\Om)} \sum_{z \in \mathbb{Z}^{n}_{\varepsilon, \Om'}}\eps^{n} \int_{Y} | \nabla m_{j} (y) |^{2}\, \di y \geq \frac{\mathcal{L}^{n}(\Om')}{\mathcal{L}^{n}(\Om)} \int_{Y}  | \nabla m_{j} (y) |^{2}\, \di y\,. \nonumber
\end{align}
Thus, $m_{j}$ is bounded in $H^{1}_{\#}(Y)$ and, up to a further subsequence, we may assume $m_{j} \rightharpoonup m$ weakly in $H^{1}_{\#}(Y)$. Clearly, $1 \leq m \leq 2 $ in $Y$ and $\fint_{\Om} m (\frac{x}{\eps}) \, \di x \leq W$. Moreover, since $\nabla (m_{j} (\frac{\cdot}{\eps})) = \frac{1}{\eps} \nabla{m}_{j} (\frac{\cdot}{\eps})$, we also have that $\nabla m_{j} (\frac{\cdot}{\eps}) \rightharpoonup \nabla m (\frac{\cdot}{\eps})$ weakly in $L^{2}(\Om; \M^{n})$. As $\C(\varphi_{j}, m_{j}) (x, \frac{x}{\eps})$ converges pointwise to $\C(\varphi, m) (x, \frac{x}{\eps})$ and is bounded in $L^{\infty}(\Om)$, we deduce that $u$ is the minimizer of~$\E_{\eps}(\varphi, m, \cdot)$ in $H^{1}_{\Gamma_{D}} (\Om; \R^{n})$. By lower-semicontinuity of~$\J_{\eps}$, $(\varphi, m, u)$ is a solution to~\eqref{e:pb-eps}.
\end{proof}

\section{The homogenized problem}
\label{s:hom}

We now introduce the homogenized versions of~$\E_{\varepsilon}$ and of $\J_{\eps}$. To do this, for $(\varphi, m) \in L^{2}(\Om; [0, 1]) \times L^{2}_{\#} (Y; [1, 2])$ we first define, for a.e.~$x \in \Om$, the \emph{correctors} (see, e.g.,~\cite{Allaire-92}) $\{w_{ij}\}_{i, j=1}^{n} \subseteq  H^{1}_{\#}(Y; \R^{n})$ relative to~$(\varphi, m)$ as the solutions to the minimum problem
\begin{align}
\label{e:min-corrector}
\min\, \Big\{ \frac12  \int_{Y}  & \C(\varphi, m)  ( \strain_{y}   (w)    +
                                   e_{ij}) {\, \cdot\,} ( \strain_{y}
                                    (w)   + e_{ij}) \, \di y : \,  w \in H^{1}_{\#} (Y; \R^{n}), \, \int_{Y} w \, \di y = 0 \Big\}, 
\end{align}
where $e_{ij} \coloneq e_{i} \odot e_{j}  =( (e_i{ \otimes }e_j) +
(e_j {\otimes }e_i))/2$,  $\{e_{i}\}_{i=1}^{n}$ being the canonical
basis of~$\R^{n}$,  and $ \strain_{y}(\cdot) $  indicates the
symmetrized gradient with respect to $y$.  

\begin{remark}
\label{r:0mean}
We notice that the correctors are well-defined in~$H^{1}_{\#} (Y; \R^{n})$ by Poincar\'e and Korn inequalities for $Y$-periodic functions with zero mean. We further point out that $w_{ij}$ also solves the minimum problem
\begin{displaymath}
\min\, \Big\{ \frac12  \int_{Y}   \C(\varphi, m)  ( \strain_{y}   (w)  +
e_{ij}) {\, \cdot\,} ( \strain_{y}  (w)  + e_{ij}) \, \di y : \,  w \in H^{1}_{\#} (Y; \R^{n})\Big\}, 
\end{displaymath}
where we drop the $0$-mean constraint. We will use this fact in Section~\ref{s:opt-sharp}.
\end{remark} 

We  now deduce a  uniform bound  on the correctors.  Since $w_{ij}$ are pointwise a.e.~defined in~$\Om$, we consider them from now on as functions $w_{ij} \colon \Om \to H^{1}_{\#} (Y; \R^{n})$.

\begin{lemma}[Bound on correctors]\label{l:correctors}
The corrector $w_{ij}$ belongs to $L^{\infty}(\Om; H^{1}_{\#}(Y; \R^{n}))$ for every $i, j=1, \ldots, n$. Moreover, 
\begin{equation}\label{e:bdd}
\sup_{\substack{ \varphi \in H^{1}(\Om; [0, 1]), \\ m \in H^{1}_{\#}(Y; [1, 2])}} \| w_{ij}\|_{L^{\infty}( \Om; H^{1}_{\#}(Y; \R^{n}))} < \infty\,.
\end{equation}
\end{lemma}

\begin{proof}
Thanks to~\eqref{e:hp-C} we deduce that for a.e.~$x \in \Om$
\begin{displaymath}
\alpha \| \strain_{y} w_{ij}(x, \cdot) + e_{ij} \|_{L^{2}_{\#}(Y ; \R^{n})}^{2} \leq 2\beta \mathcal{L}^{n}( \Om )   \qquad \text{for every $i, j=1, \ldots, n$}.
\end{displaymath}
Hence, $w_{ij}$ belongs to $L^{\infty}(\Om; H^{1}_{\#}(Y; \R^{n}))$ and \eqref{e:bdd} holds by Poincar\'e and Korn inequality.
\end{proof}

For later use, we define the homogenized tensor
\begin{align}
\label{e:C-hom}
\C^{*}_{ijkl}(\varphi, m) (x) &  \coloneq \int_{Y} \C(\varphi, m) ( \strain_{y} (w_{ij})  + e_{ij}) {\, \cdot\,} (\strain_{y} (w_{kl})  + e_{kl}) \, \di y 
\\
&
=  \int_{Y} \C(\varphi, m) ( \strain_{y} (w_{ij})  + e_{ij}) {\, \cdot\,}  e_{kl} \, \di y \,, \nonumber
\end{align}
the last equality being a consequence of~\eqref{e:min-corrector}. In particular, $\C^{*}(\varphi, m)$ is still positive definite and bounded independently of~$\varphi \in L^{2}(\Om; [0, 1])$ and of~$m \in L^{2}_{\#} (Y; [1, 2])$ (see, e.g.,~\cite[Theorem~1.4.2]{MR1859696}): for every $\varphi \in L^{2}(\Om; [0, 1])$, every~$m \in L^{2}_{\#} (Y; [1, 2])$, and a.e.~$x \in \Om$ it holds
\begin{equation}\label{e:C*}
\alpha | \e|^{2} \leq \C^{*}(\varphi, m) (x) \e {\, \cdot\,} \e \leq \beta | \e|^{2} \qquad \text{for every $\e \in \M^{n}_{S}$}.
\end{equation}

We now show the continuity of the correctors~$w_{ij}$ w.r.t.~the data~$(\varphi, m)$.

\begin{lemma}[Correctors are Lipschitz]
\label{l:cont-correctors}
The correctors $w_{ij}$ defined in~\eqref{e:min-corrector} are Lipschitz continuous in $L^{2}(\Om; H^{1}_{\#}(Y; \R^{n}))$ as functions of~$(\varphi, m) \in L^{2}(\Om; [0, 1]) \times L^{2}_{\#}(Y; [1, 2])$.
\end{lemma}

\begin{proof}
Let us denote by~$w^{\varphi, m}_{ij}$ and by $w^{\psi, \mu}_{ij}$ the correctors relative to~$(\varphi, m)$ and~$(\psi, \mu)$, respectively. By~\eqref{e:hp-C}, for a.e.~$x \in \Om$ we have that
\begin{align*}
\alpha \| \strain_{y} & ( w^{\varphi, m}_{ij}) (x, \cdot) -  \strain_{y} ( w^{\psi, \mu}_{ij}) (x, \cdot) \|_{L^{2}_{\#}(Y)}^{2} 
\\
&
\leq \int_{Y} \C(\varphi, m) (\strain_{y}  (w^{\varphi, m}_{ij})
    (x, y) -  \strain_{y} ( w^{\psi, \mu}_{ij}) (x, y) ) {\,
  \cdot\,}(\strain_{y} ( w^{\varphi, m}_{ij} )  (x, y) -  \strain_{y} ( w^{\psi, \mu}_{ij}) (x, y) )\, \di y \nonumber
\\
&
= \int_{Y} \big( \C(\psi, \mu) - \C(\varphi, m) \big) \strain_{y} (
  w^{\psi, \mu}_{ij}) (x, y) {\, \cdot\,} (\strain_{y} 
  (w^{\varphi, m}_{ij})  (x, y) -  \strain_{y} ( w^{\psi, \mu}_{ij}) (x, y) )\, \di y \nonumber
\\
&
\leq  ( 3 \beta |\varphi (x)  - \psi (x) | + \beta \|m - \mu \|_{L^{2}_{\#}(Y)})  \| \strain_{y}  ( w^{\varphi, m}_{ij}) (x, \cdot) -  \strain_{y} ( w^{\psi, \mu}_{ij}) (x, \cdot) \|_{L^{2}_{\#}(Y)} \,.\nonumber
\end{align*}
 We conclude by integrating on~$\Om$.
\end{proof}

The homogenized energy reads as
\begin{align}
\E ( \varphi, m,  u)  & \coloneq  \frac12\int_{\Om}\int_{Y} \C(\varphi, m) (\strain_{x}( u) + (\mathbf{E}_{y} w)^{\top} \strain_{x}( u) {\, \cdot\,} ( \strain_{x}( u) + (\mathbf{E}_{y} w)^{\top} \strain_{x}( u)) \, \di x 
-   \mathcal{C}(\varphi,u) 
\nonumber
\\
&
= \frac12 \int_{\Om} \C^{*}(\varphi, m) \strain_{x}(u) {\, \cdot\,}
  \strain_{x} (u) \, \di x -   \mathcal{C}(\varphi,u) 
,  \label{e:energy-2}
\end{align}
where $w = ( w_{ij})_{i,j = 1, \ldots, n} \in L^{\infty}(\Om;
H^{1}_{\#} ( Y; \M^{n}))$  are the above-introduced correctors  and $\mathbf{E}_{y} w$ is the 4-th order tensor defined componentwise by
\begin{displaymath}
(\mathbf{E}_{y}w)_{ijkl} \coloneq (\strain_{y} (w_{ij}))_{kl}\,.
\end{displaymath}
Correspondingly, the set of admissible phase fields and micro-structures writes as
\begin{displaymath}
 \A \coloneq \left\{ (\varphi, m) \in H^{1}(\Om; [0,1]) \times H^{1}_{\#}(Y; [1, 2]): \, \int_{\Om} \varphi\, \di x \leq V \text{ and } \int_{Y} m \, \di y \leq W \right\}
\end{displaymath}
and the homogenized cost functional is
\begin{equation}
\label{e:cost}
\J( \varphi, m, u) \coloneq   \mathcal{C}(\varphi,u)
+ \frac12 \int_{Y} | \nabla m|^{2} \, \di y + \frac12 \int_{\Om} | \nabla \varphi|^{2} \, \di x 
\end{equation}
for $u \in H^{1}_{\Gamma_{D}} (\Om; \R^{n})$ and $(\varphi, m) \in
 \A$,   while we set $\J(\varphi, m, u) \coloneq 
\infty$ for $(\varphi, m)  \in ( L^{2}(\Om; [0,1]) \times
L^{2}_{\#}(Y; [1, 2]) ) \setminus  \A$. 


The following proposition justifies the definition of~$\E$ as the homogenization of~$\E_{\varepsilon}$.

\begin{proposition}[$\Gamma$-convergence of $\E_{\varepsilon}$]
\label{p:Gamma-E}
 The  sequence of functionals~$\E_{\varepsilon}$
in~\eqref{e:energy} $\Gamma$-converges to~$\E$ in~\eqref{e:energy-2}
 w.r.t. the strong$\times$strong$\times$weak topology of
$L^{2}(\Om; [0,1]) \times L^{2}_{\#}(Y; [1, 2]) \times H^{1}_{\Gamma_{D}} (\Om ; \R^{n} )$.
\end{proposition}

\begin{proof}
The  assertion  
follows by a standard result  on  $\Gamma$-convergence of quadratic forms on~$H^{1}_{\Gamma_{D}} (\Om)$. Indeed, we have that for $(\varphi_{\eps}, m_{\eps})  \to (\varphi, m)$ in $L^{2}(\Om; [0,1]) \times L^{2}_{\#}(Y; [1, 2])$ the sequence of functionals
\begin{align*}
\mathcal{F}_{\varepsilon} (u) := \mathcal{E}_{\eps} (\varphi_{\eps}, m_{\eps}, u) \qquad \text{$u \in H^{1}_{\Gamma_{D}} (\Om; \R^{n})$}
\end{align*} 
$\Gamma$-converges w.r.t.  the weak topology of  $H^{1}_{\Gamma_{D}} (\Om; \R^{n})$ to
\begin{align*}
\mathcal{F} (u) := \mathcal{E} (\varphi, m, u) \qquad \text{$u \in H^{1}_{\Gamma_{D}} (\Om; \R^{n})$.}
\end{align*}
This can be deduced by~\cite[Chapter~12]{MR1201152} and~\cite[Section~2]{Allaire-92} (see also~\cite[Section~1.4]{MR1859696}), upon noticing that 
\begin{equation}
\label{e:C1}
\lim_{\eps \to 0} \int_{\Om} \Big| \C(\varphi_{\eps}, m_{\eps}) \Big( x, \frac{x}{\eps} \Big) \Big|^{2} \di x = \int_{\Om} \int_{Y} | \C(\varphi, m) (x, y) |^{2} \di y \, \di x \,.
\end{equation}
In order to prove~\eqref{e:C1}, we first show that
\begin{equation}\label{e:10}
\lim_{\eps \to 0} \int_{\Om} \Big| \C(\varphi , m_{\eps} ) \Big( x, \frac{x}{\eps} \Big) \Big|^{2} \di x = \int_{\Om} \int_{Y} | \C(\varphi, m) (x, y) |^{2} \di y \, \di x \,.
\end{equation}
With the same notation of~\eqref{e:12}, for every $\Om' \Subset \Om$  and for $\eps>0$ small enough we have that $Q^{\eps}_{z} \subseteq \Om$ for every $z \in \mathbb{Z}^{n}_{\varepsilon, \Om'}$. Hence, for such $\eps$ it holds
\begin{align}\label{e:11}
\int_{\Om}& \  \Big| \C(\varphi , m_{\eps} ) \Big( x, \frac{x}{\eps} \Big) \Big|^{2} \di x \geq \sum_{z \in \mathbb{Z}^{n}_{\varepsilon, \Om'}} \int_{Q^{\eps}_{z}}  \Big| \C(\varphi , m_{\eps} ) \Big( x, \frac{x}{\eps} \Big) \Big|^{2} \di x
\\
&
 =  \sum_{z \in \mathbb{Z}^{n}_{\varepsilon, \Om'}} \eps^{n} \int_{Y}  | \C(\varphi , m_{\eps} ) ( \eps z + \eps y, z + y ) |^{2} \di y 
=  \sum_{z \in \mathbb{Z}^{n}_{\varepsilon, \Om'}} \! \eps^{n} \! \int_{Y}  | \C(\varphi , m_{\eps} ) ( \eps z + \eps y,  y ) |^{2} \di y \nonumber
\\
&
=  \sum_{z \in \mathbb{Z}^{n}_{\varepsilon, \Om'}} \int_{Q^{\eps}_{z}} \int_{Y}  \Big| \C(\varphi , m_{\eps} ) \Big( \eps \left\lfloor \frac{x}{\eps} \right\rfloor + \eps y,  y \Big) \Big|^{2} \di y \, \di x  \nonumber
\\
&
\geq \int_{\Om'} \int_{Y} \big| T_{\eps}\varphi (x, y) m_{\eps}(y)\C_{1} + (1 - T_{\eps} \varphi (x, y)) \C_{2} \big|^{2}\, \di y \, \di x\,, \nonumber
\end{align}
where we have introduced the unfolding operator $T_{\eps} \varphi(x, y) = \varphi ( \eps \lfloor \frac{x}{\eps} \rfloor + \eps y)$ for $(x, y) \in \Om \times Y$ (see~\cite{Fon-Kro2010}). Since $\varphi \in L^{\infty}(\Om;[0,1])$, then $T_{\eps} \varphi \to \varphi$ in $L^{2}(\Om \times Y)$ (see~\cite[Proposition~A.1]{Fon-Kro2010}). Hence, passing to the liminf in~\eqref{e:11} as $\eps \to 0$ we deduce, by the arbitrariness of~$\Om' \Subset \Om$, that
\begin{displaymath}
\liminf_{\varepsilon \to 0} \int_{\Om} \Big| \C(\varphi , m_{\eps} ) \Big( x, \frac{x}{\eps} \Big) \Big|^{2} \di x \geq \int_{\Om} \int_{Y} | \C(\varphi, m) (x, y) |^{2} \di y \, \di x \,.
\end{displaymath}
The opposite inequality for the limsup can be obtained in a similar way, since the sequence $\C(\varphi , m_{\eps} ) ( x, \frac{x}{\eps}) $ is bounded in $L^{\infty}(\Om)$. Thus, we conclude for~\eqref{e:10}. We further estimate
\begin{align*}
\int_{\Om} & \  \Big| \C(\varphi_{\eps}, m_{\eps}) \Big(x, \frac{x}{\eps} \Big) -   \C(\varphi , m_{\eps}) \Big(x, \frac{x}{\eps} \Big)\Big|^{2}\di x 
\\
&
=  \int_{\Om}  \Big| (\varphi_{\eps} (x) - \varphi(x))  \Big( m_{\eps} \Big(\frac{x}{\eps} \Big) \C_{1} - \C_{2} \Big) \Big|^{2} \di x 
\leq C \| \varphi_{\eps} - \varphi\|_{L^{2}(\Om)}^{ 2}
\end{align*}
for some positive constant~$C$ independent of~$\varepsilon$. Hence,~\eqref{e:C1} holds and the two-scale convergence arguments of~\cite[Section~2]{Allaire-92} can be used to deduce that the $\Gamma$-limit of~$\mathcal{F}_{\eps}$ is the functional~$\mathcal{F}$.

We now conclude for the $\Gamma$-convergence of~$\E_{\eps}$ to~$\E$. The $\Gamma$-liminf inequality immediately follows from the~$\Gamma$-convergence of~$\F_{\eps}$ to~$\F$. As for the $\Gamma$-limsup inequality, we can fix $\varphi_{\eps} := \varphi$, $m_{\eps}:= m$, and define~$u_{\eps}$ as the recovery sequence for~$\mathcal{F}_{\eps}$,~$\mathcal{F}$. 
\end{proof}

Let us denote by~$S \colon L^{2}(\Om; [0, 1]) \times  L^{2}_{\#}(Y; [ 1, 2]) \to H^{1}_{\Gamma_{D}}(\Om; \R^{n})$ the map that to each $(\varphi, m) \in L^{2}(\Om; [0, 1]) \times L^{2}_{\#}(Y; [ 1, 2])$ associates the unique solution to
\begin{equation}
\label{e:S}
\min_{u \in H^{1}_{\Gamma_{D}} (\Om; \R^{n})} \E(\varphi, m, u)\,.
\end{equation}
To shorten the notation, let us set 
$$\mathcal{G}_{\eps} (\varphi, m)
\coloneq \mathcal{J}_{\eps} (\varphi, m, S_{\eps}(\varphi, m)) \ \
\text{and} \ \ \mathcal{G} (\varphi, m) \coloneq \mathcal{J} (\varphi,
m, S(\varphi, m))$$ 
for $(\varphi, m) \in H^{1}(\Om; [0, 1]) \times H^{1}_{\#}(Y; [ 1,
2])$. We recall that, by the definition of~$\mathcal{J}$ and
of~$\mathcal{J}_{\varepsilon}$, it holds $\G_{\varepsilon} (\varphi,
m) =  \infty$ for $(\varphi, m) \in ( L^{2}(\Om; [0, 1])
\times L^{2}_{\#} (Y; [1, 2]) ) \setminus 
\mathcal{A}_{\varepsilon}$  and $\G (\varphi, m) = 
\infty$ for every $(\varphi, m) \in (  L^{2}(\Om; [0, 1]) \times
L^{2}_{\#}(Y; [ 1, 2])) \setminus  \A$, Then, the following holds.

\begin{proposition}[$\Gamma$-convergence of $\mathcal{G}_{\eps}$]
\label{p:Gamma-J}
 The sequence of functionals  $\mathcal{G}_{\eps}$
$\Gamma$-converges to~$\mathcal{G}$  with respect to the
strong  topology of  $L^{2}(\Om; [0,1]) \times
L^{2}_{\#}(Y;[1,2])$. 
\end{proposition}

\begin{proof}
Let $\varphi_{\eps} \to \varphi$ in $L^{2}(\Om;[0,1])$ and $m_{\eps}
\to m$ in $L^{2}_{\#}(Y;[1,2])$ be  given in such  a way that 
$\liminf_{\eps \to 0} \mathcal{G}_{\eps} (\varphi_{\eps}, m_{\eps})
< \infty$. Up to a subsequence, we may assume that the liminf
is a limit, that $\mathcal{G}_{\eps}(\varphi_{\eps}, m_{\eps})$ is
bounded, and that $(\varphi_{\eps}, m_{\eps}) \in  \A_{\eps}$.  Hence, $\varphi_{\eps} \rightharpoonup \varphi$ weakly in $H^{1}(\Om)$ and $\int_{\Om} \varphi \, \di x \leq V$. For every $\Om' \Subset \Om$, arguing as in~\eqref{e:12} and in~\eqref{e:11} we get
\begin{align}\label{e:13}
\fint_{\Om} \Big| \nabla m_{\eps}\Big(\frac{x}{\eps} \Big) \Big|^{2}\, \di x & \geq \frac{\mathcal{L}^{n}(\Om')}{\mathcal{L}^{n}(\Om)} \int_{Y}  | \nabla m_{\eps} (y) |^{2}\, \di y\,. 
\end{align}
Thus, $m_{\eps} \rightharpoonup m$ weakly in~$H^{1}_{\#}(Y)$ and, being~$\Om'\Subset \Om$ arbitrary in~\eqref{e:13}, we have that
\begin{displaymath}
\liminf_{\eps \to 0} \fint_{\Om}\Big| \nabla m_{\eps}\Big(\frac{x}{\eps} \Big) \Big|^{2}\, \di x \geq \int_{Y}  | \nabla m (y) |^{2}\, \di y\,.
\end{displaymath}
In the same way, for every $\Om'' \Supset \Om$ we can show that
\begin{displaymath}
\fint_{\Om}  m_{\eps}\Big(\frac{x}{\eps} \Big) \, \di x \leq \frac{\mathcal{L}^{n}(\Om'')}{\mathcal{L}^{n}(\Om)} \int_{Y}   m_{\eps} (y) \, \di y\,,
\end{displaymath}
 which  in turn implies that $\int_{Y} m \, \di y \leq
W$. Thus, $(\varphi, m) \in \A$.  By Proposition~\ref{p:Gamma-E} we have that $S_{\eps} (\varphi_{\eps}, m_{\eps}) \rightharpoonup S(\varphi, m)$ weakly in~$H^{1}_{\Gamma_{D}} (\Om; \R^{n})$. Hence, 
\begin{displaymath}
\mathcal{G}(\varphi, m) \leq \liminf_{\eps\to 0} \mathcal{G}_{\eps} (\varphi_{\eps}, m_{\eps})\,.
\end{displaymath}
A recovery sequence can be simply taken to be constant.
\end{proof}

In view of Propositions~\ref{p:Gamma-E} and~\ref{p:Gamma-J}, in what
follows we consider the {\it Homogenized  Microstructure-Topology
  Optimization} (HMTO)  problem 

\begin{align}
\label{e:pb}
& \min_{\substack{\varphi \in H^{1}(\Om; [0, 1]),\\ m \in
  H^{1}_{\#}(Y; [ 1, 2]) }} \, \J(\varphi, m, u)  \ \ \  \text{subject
  to } \min_{u \in H^{1}_{\Gamma_{D}} (\Om; \R^{n}) } \, \E(\varphi,
  m, u) \,. \tag{HMTO}
\end{align}

\begin{proposition}[Existence for the (HMTO) problem]
\label{p:existence}
Problem~\eqref{e:pb} admits a solution. Moreover, for every
sequence~$(\varphi_{\eps}, m_{\eps}, u_{\eps})$ of  solutions  to~\eqref{e:pb-eps} there exists a solution $(\varphi, m, u)$ to~\eqref{e:pb} such that, up to a subsequence, $\varphi_{\eps} \rightharpoonup \varphi$ weakly in $H^{1}(\Om)$, $m_{\eps} \rightharpoonup m$ weakly in $H^{1}_{\#}(Y)$, and $u_{\eps} \to u$ in~$H^{1}_{\Gamma_{D}} (\Om; \R^{n})$.
\end{proposition}

\begin{proof}
The existence of a minimizer of~\eqref{e:pb} can be proved by the
Direct Method: Let $(\varphi_{j}, m_{j}, u_{j}) \in H^{1}(\Om; [0,1])
\times H^{1}_{\#}(Y; [1, 2]) \times H^{1}_{\Gamma_{D}} (\Om; \R^{n})$
be a minimizing sequence for~\eqref{e:pb}. Then, up to a  not
relabeled  subsequence, we may assume that $(\varphi_{j}, m_{j})
\in  \A$,   $\varphi_{j} \rightharpoonup \varphi$ weakly
in~$H^{1}(\Om)$, and $m_{j} \rightharpoonup m$ weakly
in~$H^{1}_{\#}(Y)$. Clearly, $0 \leq \varphi \leq 1$, $1 \leq m \leq
2$,  so that  $(\varphi, m) \in \A$.   Lemmas~\ref{l:correctors}-\ref{l:cont-correctors} imply that, up to a not relabeled subsequence, $\C^{*}(\varphi_{j}, m_{j})$ converges to~$\C^{*}(\varphi, m)$ pointwise a.e.~in~$\Om$ and is bounded in $L^{\infty}(\Om)$, which implies that $u_{j} \to u:= S(\varphi, m)$ in~$H^{1}(\Om; \R^{n})$. Finally, the lower semicontinuity of~$\J$ implies that~$(\varphi, m, u)$ is a solution to~\eqref{e:pb}.

The second part of the proposition is a consequence of the $\Gamma$-convergence results of Propositions~\ref{p:Gamma-E}-\ref{p:Gamma-J} and standard arguments in linear elasticity.
\end{proof}

\subsection{Optimality conditions}
\label{s:optimality}

 In this section, we prove some  first-order optimality
conditions  for the \eqref{e:pb} problem.  First, let us
define  a suitable  class of {\it admissible variations}.

\begin{definition}[Admissible variations]\label{d:a}
We say that the pair $(\psi, \mu) \in H^{1}(\Om) \times H^{1}_{\#}
(Y)$ is an \emph{admissible variation} for~$(\varphi, m)$, and we
write $(\psi, \mu) \in \V(\varphi, m)$, if there exists
$\overline{t}>0$ such that for every $t \in [0, \overline{t}]$ the
pair $(\varphi +  t\psi, m+t\mu)$ belongs to $\A$.  
\end{definition}
 
Notice that, since $0 \leq \varphi \leq 1$ and $1\leq m \leq 2$, every
admissible variation $(\psi, \mu) \in {\rm V} (\varphi, m)$ belongs to
$L^{\infty}(\Om) \times L^{\infty}_{\#}(Y)$. 

The main result of the section reads as follows.

\begin{theorem}[Optimality conditions for~\eqref{e:pb}]
\label{t:optimality-deterministic}
Let 
$(\varphi, m, u) \in H^{1}(\Om; [0,1]) \times H^{1}_{\#}(Y;[1,2]) \times H^{1}_{\Gamma_{D}}(\Om; \R^{n})$ be a solution to~\eqref{e:pb}. Then, for every $v \in H^{1}_{\Gamma_{D}} (\Om; \R^{n})$ and every $(\psi, \mu) \in \V(\varphi, m)$ the following optimality conditions hold:
\begin{align}
& \int_{\Om} \C^{*}(\varphi, m) \strain_{x} (u) {\, \cdot\,}
  \strain_{x}(v) \, \di x =   \mathcal{C}(\varphi,u)
\,, \label{e:Copt-1}
\\
&
2 \int_{\Om} \psi f {\, \cdot\,} u \, \di x - \int_{\Om} \overline{\C}^{*}(\psi, \mu) \strain_{x} (u) {\, \cdot\,} \strain_{x}(u) \, \di x + \int_{\Om} \nabla \varphi{\, \cdot\,} \nabla \psi \, \di x + \int_{Y} \nabla m {\, \cdot\,} \nabla \mu \, \di y \geq 0\,, \label{e:opt-2}
\end{align}
where for $i, j, k, l = 1, \ldots, n$ we have set
\begin{equation}
\label{e:Cbar}
\overline{\C}^{*}_{ijkl}(\psi, \mu) \coloneq \int_{Y} \big( \psi(m \C_{1} - \C_{2} ) + \varphi \mu \C_{1} \big) (\strain_{y}(w_{ij} )   + e_{ij}) {\, \cdot\,} (\strain_{y}(w_{kl})  + e_{kl})\, \di y\,.
\end{equation}
\end{theorem}

 The first step towards the proof of Theorem~\ref{t:optimality-deterministic} is to show the differentiability of the homogenized elasticity tensor~$\C^{*}$ w.r.t.~variations in $\mathrm{V}(\varphi, m)$. This is actually a consequence of the following proposition.

\begin{proposition}[Differentiability of $w^t_{ij}$]
\label{p:C*-diff}
Let 
$(\psi, \mu) \in \V(\varphi, m)$. Let us denote by $w^{t}_{ij} \in L^{\infty}(\Om; H^{1}_{\#} (Y; \R^{n}))$ the correctors relative to $(\varphi + t\psi, m + t \mu)$. Then, the function $t\mapsto w^{t}_{ij}$ is differentiable in $t=0$ and the derivative $\overline{w}_{ij}(\psi, \mu) = \overline{w}_{ij} \coloneq \partial_{t}|_{t=0} w^{t}_{ij} \in L^{\infty}(\Om; H^{1}_{\#}(Y; \R^{n}))$ satisfies for a.e.~$x \in \Om$
\begin{align}
\label{e:overlinew}
\int_{Y} \C(\varphi, m) & \strain_{y} (\overline{w}_{ij}) {\, \cdot\,} \strain_{y} (v) \, \di y 
\\
&
+  \int_{Y}\big[\psi ( m\C_{1} - \C_{2}) + \varphi \mu \C_{1} \big] ( \strain_{y} (w_{ij}) + e_{ij}) {\, \cdot\,} \strain_{y}(v) \, \di y = 0 \nonumber
\end{align}
for every $v \in  H^{1}_{\#}(Y; \R^{n})$ with $\int_{Y} v \, \di y = 0$.
\end{proposition}

\begin{remark}
 The Lax-Milgram lemma ensure that  solution
$\overline{w}_{ij}$ to~\eqref{e:overlinew}  uniquely  exists in the space of functions $w \in  H^{1}_{\#}(Y; \R^{n})$ such that $\int_{Y} w \, \di y = 0$. Interpreting $\overline{w}_{ij}$ as a function from $\Om$ to $H^{1}_{\#}(Y; \R^{n})$, it is possible to show that $\overline{w}_{ij} \in L^{\infty} (\Om; H^{1}_{\#}(Y; \R^{n}))$.
\end{remark} 

\begin{proof}
By Lemma~\ref{l:cont-correctors}, since $(\psi, \mu) \in L^{\infty}(\Om) \times L^{\infty}_{\#}(Y)$ we obtain
\begin{equation}\label{e:C3}
\| w_{ij} - w^{t}_{ij} \|_{L^{\infty}(\Om; H^{1}_{\#}(Y; \R^{n}))} \leq C t (\| \psi\|_{L^{\infty}(\Om)} + \| \mu\|_{L^{\infty}_{\#}(Y)})
\end{equation}
for some positive constant~$C$ independent of~$t$, $\psi$, and~$\mu$.

We now show that 
\begin{displaymath}
\| w^{t}_{ij} - w_{ij} - t \overline{w}_{ij}\|_{L^{\infty}(\Om;
  H^{1}_{\#}(Y; \R^{n}))} \leq  {\rm O}(t^2)  \,,
\end{displaymath}
which indeed implies the differentiability of~$t\mapsto w^{t}_{ij}$ in
$t=0$,  as well as the fact that $\overline{w}_{ij}
= \partial_t|_{t=0}w^t_{ij}$.  For a.e.~$x \in \Om$ we have that
\begin{align}
\label{e:C4}
0 = & \int_{Y} \C(\varphi + t\psi, m + t\mu) (\strain_{y} (w^{t}_{ij}) + e_{ij} ){\, \cdot\,} \strain_{y} (v)  \, \di y 
\\
&
 - \int_{Y} \C(\varphi , m ) (\strain_{y} (w_{ij}) + e_{ij}) {\, \cdot\,} \strain_{y} (v)  \, \di y  
- t \int_{Y} \C(\varphi, m)  \strain_{y} (\overline{w}_{ij}) {\, \cdot\,} \strain_{y} (v) \, \di y \nonumber
\\
&
-  t  \int_{Y}\big[\psi ( m \C_{1} - \C_{2}) + \varphi \mu \C_{1} \big] ( \strain_{y} (w_{ij}) + e_{ij}) {\, \cdot\,} \strain_{y}(v) \, \di y \nonumber
\end{align}
for every $v \in H^{1}_{\#}(Y; \R^{n})$ with $\int_{Y} v \, \di y = 0$. Let us set $v_{t} \coloneq w^{t}_{ij} - w_{ij} - t\overline{w}_{ij}$. By simple algebraic manipulation, we rewrite~\eqref{e:C4} as
\begin{align}
\label{e:5}
0 = &  \int_{Y} \C(\varphi , m ) \strain_{y}(v_{t})  {\, \cdot\,} \strain_{y}(v) \, \di y 
\\
&
+ \int_{Y} \big( \C(\varphi + t\psi, m + t\mu) - \C(\varphi, m) - t\big[\psi ( m \C_{1} - \C_{2}) + \varphi \mu \C_{1} \big] \big) (\strain_{y} (w^{t}_{ij}) + e_{ij} ){\, \cdot\,} \strain_{y} (v)  \, \di y \nonumber
\\
&
+ t\int_{Y} \big[\psi ( m\C_{1} - \C_{2}) + \varphi \mu \C_{1} \big] ( \strain_{y} (w^{t}_{ij}) - \strain_{y} (w_{ij})) {\, \cdot\,} \strain_{y} (v) \, \di y \,.\nonumber
\end{align}
Choosing $v = v_{t}$ in \eqref{e:5}, from~\eqref{e:C3} we deduce that for a.e.~$x \in \Om$
\begin{align*}
\alpha \| \strain_{y} (v_{t}) \|_{L^{2}_{\#}(Y)} \leq Ct^{2} \big(\| \varphi\|_{L^{\infty}(\Om)} \| \mu\|_{L^{\infty}_{\#}(Y)} + \| \psi\|_{L^{\infty}(\Om)} + \| \mu\|_{L^{\infty}_{\#}(Y)} \big) 
\end{align*}
for a positive constant~$C$ independent of~$t$, $\mu$,
and~$\psi$. Since $\int_{Y} v_{t} \, \di y = 0$ a.e.~in~$\Om$, by
Korn's inequality in $H^{1}_{\#}(Y; \R^{n})$ and taking the supremum
over $\Om$ we infer that $\| v_{t}\|_{L^{\infty}(\Om; H^{1}_{\#}(Y;
  \R^{n}))} =  {\rm O}(t^{2})$ and $t \mapsto w^{t}_{ij}$ is differentiable in $t=0$ for every $i, j=1, \ldots, n$.
\end{proof}

We now conclude the differentiability of the homogenized elasticity tensor~$\C^{*}$.

\begin{corollary}[Differentiability of $\C^{*}$]\label{c:C*-diff}
Under the assumptions of Proposition~\emph{\ref{p:C*-diff}}, the function~$t \mapsto \C^{*} (\varphi + t\psi, m+ t\mu)$ is differentiable in $t=0$ in $L^{\infty}(\Om)$ and, for every $i, j, k, l= 1, \ldots, n$,
\begin{align}\label{e:C*-diff}
 \partial_{t}|_{t=0} \C^{*} (\varphi + t\psi, m+t\mu) = \overline{\C}^{*}(\psi, \mu) \,,
 \end{align}
 where $\overline{\C}^{*}(\psi, \mu)$ is defined in~\eqref{e:Cbar}.
\end{corollary}

\begin{proof}
Let us fix $(\varphi, m) \in  \A$,  $(\psi, \mu) \in \V(\varphi,m)$, and $t>0$ small enough, and let $\overline{w}_{ij} \in L^{\infty}(\Om; H^{1}_{\#}(Y; \R^{n}))$ be the functions identified in Proposition~\ref{p:C*-diff}. Recalling the definition~\eqref{e:min-corrector} of the correctors $w_{ij}$, for a.e.~$x \in \Om$ we have that
\begin{align*}
& \int_{Y} \C(\varphi, m) \strain_{y} ( \overline{w}_{ij} ) {\, \cdot\,} (\strain_{y}(w_{kl}) + e_{kl}) \, \di y = 0\,,
\\
&  \int_{Y} \C(\varphi, m)  (\strain_{y}(w_{ij}) + e_{ij}){\, \cdot\,}  \strain_{y} (\overline{w}_{kl})\, \di y = 0\,.
\end{align*} 
Thus, we estimate pointwise  almost everywhere  in $\Om$ the difference
\begin{align*}
\C^{*}_{ijkl} & (\varphi + t\psi, m+t\mu) - \C^{*}_{ijkl}(\varphi, m) - t\overline{\C}^{*}_{ijkl}(\psi, \mu)  
\\
&
= \int_{Y} \C(\varphi + t\psi, m + t\mu) ( \strain_{y} (w^{t}_{ij})  + e_{ij}) {\, \cdot\,} (\strain_{y} (w^{t}_{kl})  + e_{kl}) \, \di y
\\
&
 \quad - \int_{Y} \C(\varphi, m) ( \strain_{y} (w_{ij})  + e_{ij}) {\, \cdot\,} (\strain_{y} (w_{kl})  + e_{kl}) \, \di y 
 \\
 &
 \quad - t   \int_{Y} \C(\varphi, m) \strain_{y} ( \overline{w}_{ij} ) {\, \cdot\,} (\strain_{y}(w_{kl}) + e_{kl}) \, \di y 
 \\
 &
 \quad - t \int_{Y} \C(\varphi, m)  (\strain_{y}(w_{ij}) + e_{ij}){\, \cdot\,}  \strain_{y} (\overline{w}_{kl})\, \di y
\\
&
\quad - t \int_{Y} \big( \psi(m \C_{1} - \C_{2} ) + \varphi \mu \C_{1} \big) (\strain_{y}(w_{ij} ) + e_{ij}) {\, \cdot\,} (\strain_{y}(w_{kl}) + e_{kl})\, \di y
\\
&
= \int_{Y} \big( \C(\varphi + t\psi, m + t\mu) - \C(\varphi, m) - t
  \big( \psi(m \C_{1} - \C_{2} ) + \varphi \mu \C_{1} \big)\big)   (
  \strain_{y} (w^{t}_{ij})  + e_{ij}) {\, \cdot\,}\\
&\qquad  \qquad  {\, \cdot\,} (\strain_{y} (w^{t}_{kl})  + e_{kl}) \, \di y
\\
&
\quad + \int_{Y} \C(\varphi, m) ( \strain_{y} (w^{t}_{ij})  + e_{ij}) {\, \cdot\,} ( \strain_{y} (w^{t}_{kl})  - \strain_{y} (w_{kl}) - t\strain_{y} (\overline{w}_{kl})) \,\di y
\\
&
\quad + \int_{Y} \C(\varphi, m) (\strain_{y} (w^{t}_{ij})  - \strain_{y} (w_{ij}) - t\strain_{y} (\overline{w}_{ij})) {\, \cdot\,}  ( \strain_{y} (w_{kl})  + e_{kl}) \,\di y
\\
&
\quad +t \int_{Y} \C(\varphi, m) ( \strain_{y} (w^{t}_{ij}) - \strain_{y} (w_{ij}) ) {\, \cdot\,} \strain_{y} (\overline{w}_{kl}) \, \di y 
\\
&
\quad + t \int_{Y} \big( \psi(m \C_{1} - \C_{2} ) + \varphi \mu \C_{1} \big) ( \strain_{y} (w^{t}_{ij})  - \strain_{y} (w_{ij})) {\, \cdot\,} (\strain_{y} (w^{t}_{kl})  + e_{kl}) \, \di y
\\
&
\quad + t \int_{Y} \big( \psi(m \C_{1} - \C_{2} ) + \varphi \mu \C_{1} \big) (\strain_{y} (w_{ij})  + e_{ij}) {\, \cdot\,} ( \strain_{y} (w^{t}_{kl})  - \strain_{y} (w_{kl}))  \, \di y\,.
\end{align*}
In view of Proposition~\ref{p:C*-diff}, the previous equality implies that
\begin{align*}
\| \C^{*}_{ijkl}&   (\varphi + t\psi, m+t\mu) - \C^{*}_{ijkl}(\varphi, m) - t\overline{\C}^{*}_{ijkl} (\psi, \mu)\|_{L^{\infty}(\Om)}  \leq C t^{2}
\end{align*}
for some positive constant~$C$ independent of~$t$. Hence, $t \mapsto \C^{*}_{ijkl}( \varphi + t\psi, m+t\mu)$ is differentiable in~$t=0$ and~\eqref{e:C*-diff} holds.
\end{proof}

We now  check  the differentiability of the control-to-state
operator~$S \colon  \A   \to H^{1}_{\Gamma_{D}}(\Om; \R^{n})$
of the forward problem in~\eqref{e:pb}, that is, the map that to each
$(\varphi, m) \in  \A$ associates the unique solution~$S(\varphi, m) \in H^{1}_{\Gamma_{D}} (\Om; \R^{n})$ to
\begin{equation}\label{e:forward}
\min_{u \in H^{1}_{\Gamma_{D}} (\Om; \R^{n})} \, \E(\varphi, m, u)\,.
\end{equation}

\begin{proposition}[Differentiability of $S$]
\label{p:u-diff}
Let 
$(\psi, \mu) \in \V(\varphi, m)$. Then, the map $t \mapsto S(\varphi + t\psi, m + t\mu)$ is differentiable in~$t=0$. Moreover, the function $v_{(\psi, \mu)} \coloneq \partial_{t}|_{t=0}S(\varphi + t\psi, m + t\mu) \in H^{1}_{\Gamma_{D}}(\Om; \R^{n})$ solves
\begin{align}
\label{e:u-diff}
& \int_{\Om} \C^{*}(\varphi, m) \strain_{x} (v_{(\psi, \mu)}) {\, \cdot\,} \strain_{x}(v) \, \di x + \int_{\Om}\overline{\C}^{*}(\psi, \mu) \strain_{x} (S(\varphi, m)) {\, \cdot\,} \strain_{x} (v) \, \di x = \int_{\Om} \psi f {\, \cdot\,} v \, \di x
\end{align} 
for every $v \in H^{1}_{\Gamma_{D}}(\Om; \R^{n})$.
\end{proposition}

\begin{proof}
We notice that the variational equation~\eqref{e:u-diff} admits unique
solution by  the Lax-Milgram lemma.  We denote such solution by~$v_{(\psi, \mu)}$. For $t$ small enough, set $u_{t} \coloneq S(\varphi + t\psi, m + t\mu)$. We now proceed as in Proposition~\ref{p:C*-diff} with the estimate of the difference $v_{t} \coloneq u_{t} - u - tv_{(\psi, \mu)}$. For every $v \in H^{1}_{\Gamma_{D}} (\Om; \R^{n})$ we have
\begin{align}
\label{e:6}
0 & =  \int_{\Om} \C^{*}(\varphi + t\psi, m+t\mu) \strain_{x} (u_{t}) {\, \cdot\,} \strain_{x} (v) \, \di x - \int_{\Om} \C^{*}(\varphi, m) \strain_{x} (u) {\, \cdot\,} \strain_{x} (v) 
\\
&
\quad - t \int_{\Om} \C^{*} (\varphi, m) \strain_{x}( v_{(\psi,\mu)}) {\, \cdot\,} \strain_{x} (v) \, \di x - t\int_{\Om} \overline{\C}^{*} (\psi, \mu) \strain_{x} (u) {\, \cdot\,} \strain_{x} (v) \, \di x \nonumber
\\
&
= \int_{\Om} \C^{*}(\varphi + t\psi, m+t\mu) \strain_{x} (v_{t}) {\, \cdot\,} \strain_{x}(v) \, \di x \nonumber
\\
&
\quad + \int_{\Om} \big( \C^{*}(\varphi + t\psi, m+t\mu) - \C^{*} (\varphi, m)  - t \overline{\C}^{*} (\psi, \mu) \big) \strain_{x} (u) {\, \cdot\,} \strain_{x}(v) \, \di x \nonumber
\\
&
\quad + t \int_{\Om} \big( \C^{*}(\varphi + t\psi, m+t\mu) - \C^{*}(\varphi, m) \big) \strain_{x}( v_{(\psi, \mu)}) {\, \cdot\,} \strain_{x} (v) \, \di x \,. \nonumber
\end{align}
Testing~\eqref{e:6} with $v = v_{t}$ and applying Proposition~\ref{p:C*-diff} and Corollary~\ref{c:C*-diff} we infer that
\begin{align*}
\alpha \| v_{t} \|_{H^{1}(\Om; \R^{n})} & \leq   \| u \|_{H^{1}(\Om;
                                          \R^{n})} \|  \C^{*}(\varphi
                                          + t\psi, m+t\mu) - \C^{*}
                                          (\varphi, m)  - t
                                          \overline{\C}^{*} (\psi,
                                          \mu)
                                          \|_{L^{\infty}(\Om)} 
\\
&
 \quad + t \| v_{(\psi, \mu)} \|_{H^{1} (\Om; \R^{n})} \| \C^{*}(\varphi + t\psi, m+t\mu) - \C^{*}(\varphi, m)\|_{L^{\infty}(\Om)} 
\\
&
\leq C t^{2} ( \| u \|_{H^{1}(\Om; \R^{n})} + \| v_{(\psi, \mu)} \|_{H^{1} (\Om; \R^{n})} ) 
\end{align*}
for a positive constant~$C$ independent of~$t$. Hence, $t \mapsto u_{t}$ is differentiable in $t=0$ and $\partial_{t}|_{t=0} u_{t} = v_{(\psi, \mu)}$.
\end{proof}

We are now in a position  of concluding  the proof of Theorem~\ref{t:optimality-deterministic}.

\begin{proof}[Proof of Theorem \ref{t:optimality-deterministic}]
Equality~\eqref{e:Copt-1} is simply a rephrase of the minimization
problem~\eqref{e:forward}. Let us prove~\eqref{e:opt-2}. For $t>0$
small enough, we consider the variation $(\varphi + t\psi, m+t\mu) \in
 \A$  and define $u_{t} \coloneq S(\varphi + t\psi, m+t\mu)$ and denote by $v_{(\psi, \mu)}$ the derivative of~$u_{t}$ in $t=0$ identified in Proposition~\ref{p:u-diff}. Since~$(\varphi, m,u)$ solves~\eqref{e:pb}, we have that $\J(\varphi, m, u) \leq \J(\varphi + t\psi, m+t\mu, u_{t})$. Dividing by $t>0$ and passing to the limit as $t\searrow 0$, by Proposition~\ref{p:u-diff} we have that
\begin{align}\label{e:7}
0 \leq & \  \int_{\Om} \psi f {\, \cdot\,} u \, \di x + \int_{\Om} \varphi f {\, \cdot\,} v_{(\psi, \mu)}\, \di x + \int_{\Gamma_{N}} g{\, \cdot\,} v_{(\psi, \mu)}\, \di \HH^{n-1} 
\\
&
+ \int_{\Om} \nabla \varphi{\, \cdot\,} \nabla \psi \, \di x + \int_{Y} \nabla m {\, \cdot\,} \nabla \mu \, \di y\,. \nonumber
\end{align}
We conclude by rewriting the second and the third terms in~\eqref{e:7}. Indeed, since $v_{(\psi, \mu)} \in H^{1}_{\Gamma_{D}}(\Om; \R^{n})$, by the minimality of~$u$ and by~\eqref{e:u-diff} we have that
\begin{align}\label{e:8}
\int_{\Om} & \varphi                 f{\, \cdot\,}
                                                  v_{(\psi, \mu)}
                                                  \,\di x +
                                                  \int_{\Gamma_{N}} g
                                                  {\, \cdot\,}
                                                  v_{(\psi, \mu)}\,
                                                  \di \HH^{n-1}  =  \int_{\Om} \C^{*}(\varphi, m) \strain_{x}(u) {\, \cdot\,} \strain_{x} (v_{(\psi, \mu)}) \, \di x 
\\
&
= -  \int_{\Om}\overline{\C}^{*}(\psi, \mu) \strain_{x} (u) {\, \cdot\,} \strain_{x} (u) \, \di x + \int_{\Om} \psi f{\, \cdot\,} u \, \di x \,. \nonumber
\end{align}
Inserting~\eqref{e:8} in~\eqref{e:7} we deduce~\eqref{e:opt-2}, which concludes the proof of the theorem.
\end{proof}

\section{The sharp-interface limit}
\label{s:sharp}

This section is devoted to the study of the sharp-interface
counterpart of~\eqref{e:pb-eps}. We start by modifying the cost
functional~$\mathcal{J}_{\eps}$  by  adding a Modica-Mortola
term~\cite{MR866718, MR0445362}, which penalizes fast oscillations
of~$\varphi$ and of~$m$. We fix $H \colon \R \to [0,  \infty)$
such that $H (\zeta)=0$ if and only if~$\zeta \in \{0, 1\}$. Given 
$(\varphi, m) \in  \A $ and $u \in H^{1}_{\Gamma_{D}} (\Om; \R^{n})$ we set
\begin{align}
\label{e:Jsharp}
 \J^{ \text{s}}_\eps (\varphi, m, u)& \coloneq   
                                                 \mathcal{C}(\varphi,u)
+ \frac{\eps}{2} \int_{\Om} | \nabla \varphi|^{2} \, \di x + \frac{1}{2\eps} \fint_{\Om} | \nabla m| ^{2} \Big( \frac{x}{\eps} \Big) \, \di x 
\\
&
+ \frac{1}{2\varepsilon} \int_{\Om}   H(\varphi (x) ) \, \di x +\frac{1}{2\varepsilon} \fint_{\Om} H \Big( m\Big(\frac{x}{\eps} \Big) -1\Big) \,  \di x\,. \nonumber
\end{align} 
The functional~$\J^{ \text{s}}_\eps $ is extended to $
\infty$ for $(\varphi, m) \in ( L^{2}(\Om; [0, 1]) \times
L^{2}_{\#} (Y; [1, 2]))  \setminus \A$.   The sharp-interface limit of~$\J^{ \text{s}}_\eps $ reads as follows. We set 
\begin{displaymath}
 \A^{ \text{s}} := \left\{ (\varphi, m)   \in BV(\Om; \{0, 1\}) \times  BV_{\#} (Y; \{1, 2\}): \, \int_{\Om} \varphi\, \di x \leq V ,\, \int_{Y} m \, \di y \leq W \right\}\,.
\end{displaymath}
For $(\varphi, m) \in  \A^{ \text{s}}$  and $u \in H^{1}_{\Gamma_{D}} (\Om; \R^{n})$ we define
\begin{align}
\label{e:Jsharp-2}
\J^{ \text{s}}  (\varphi, m, u) \coloneq   \mathcal{C}(\varphi,u)
+ c_{H} \big(  \mathcal{P} (\{\varphi = 1\}; \Om) + \mathcal{P} (\{ m=2\}; Y) \big)\,,
\end{align}
where $\mathcal{P} (B; U)$ for $U$ open and $B \subseteq U$ denotes
the perimeter of~$B$ in~$U$, and~$c_{H}\coloneq \int_{0}^{1} \sqrt{2
  H(t)} \, \di t$. For $(\varphi, m) \in  ( BV(\Om; \{0, 1\}) \times
BV_{\#} (Y; \{1, 2\}) ) \setminus   \A^{ \text{s}} $  we set  $\J^{ \text{s}}  (\varphi, m, u) \coloneq  \infty$. 

Recalling the definition of~$S_{\eps}, S \colon L^{2} (\Om; [0, 1])
\times L^{2}_{\#} (Y; [1, 2]) \to H^{1}_{\Gamma_{D}} (\Om; \R^{n})$
given in~\eqref{e:pb-eps} and in~\eqref{e:S}, respectively,  we
let
\begin{align*}
\mathcal{V}_{\eps} (\varphi, m)& \coloneq \J^{ \text{s}}_\eps  (\varphi, m, S_{\eps} (\varphi, m)) \qquad \text{for $(\varphi, m) \in L^{2} (\Om; [0, 1]) \times L^{2}_{\#} (Y; [1, 2])$,}\\
\mathcal{V} (\varphi, m) & \coloneq \J^{ \text{s}} (\varphi, m, S (\varphi, m)) \qquad \text{for $(\varphi, m) \in BV(\Om; \{0, 1\}) \times BV_{\#} (Y; \{1, 2\} )$,}
\end{align*}
We  can prove the next $\Gamma$-convergence result.

\begin{proposition}[$\Gamma$-convergence of $\mathcal{V}_{\eps} $]
\label{p:sharp-conv}
 The sequence of functionals $\mathcal{V}_{\eps}$
$\Gamma$-converges to~$\mathcal{V}$  with respect to the strong
topology in $L^{1}(\Om; [0, 1]) \times L^{1}_{\#} (Y; [1, 2])$. 
\end{proposition}

\begin{proof}
By Proposition~\ref{p:Gamma-E}, we have that $S_{\eps}
(\varphi_{\eps}, m_{\eps}) \rightharpoonup S(\varphi, m)$  in
$H^1_{\Gamma_D}(\Omega;\R^n)$  whenever $(\varphi_{\eps}, m_{\eps}) \to (\varphi, m)$ in $L^{1}(\Om; [0,1]) \times L^{1}_{\#} (Y; [1, 2])$. Reasoning as in~\eqref{e:12}, we notice that for every $\Om' \Subset \Om$ and for $\eps>0$ small enough we estimate
\begin{align}
\label{e:1000C}
\frac{\eps}{2} \fint_{\Om}  | \nabla m_{\eps}| \Big(\frac{x}{\eps} \Big) \, \di x &  + \frac{1}{2\eps} \fint_{\Om} H\Big( m\Big(\frac{x}{\eps} \Big) -1\Big) \di x 
\\
&
\nonumber \geq \frac{\eps \mathcal{L}^{n} (\Om')}{2 \mathcal{L}^{n} (\Om) } \int_{Y} | \nabla m_{\eps}|^{2} (y) \, \di y + \frac{\mathcal{L}^{n} (\Om')}{2\eps \mathcal{L}^{n} (\Om)} \int_{Y} H(m_{\eps} (y) - 1) \, \di y \,.
\end{align}
Hence,~\eqref{e:1000C} and the classical Modica-Mortola~\cite{MR866718, MR0445362} argument yield the $\Gamma$-liminf inequality
\begin{displaymath}
\liminf_{\eps \to 0} \, \mathcal{V}_{\eps} (\varphi_{\eps}, m_{\eps}) \geq \mathcal{V} (\varphi, m)
\end{displaymath}
whenever $(\varphi_{\eps}, m_{\eps}) \to (\varphi, m)$ in $L^{1} (\Om;
[0, 1]) \times L^{1}_{\#} (Y; [1, 2])$ with~$(\varphi_{\eps},
m_{\eps}) \in \A$.  Arguing as in Proposition~\ref{p:Gamma-J}
(see~\eqref{e:13}) we deduce that the limit $(\varphi, m)$ belongs to
$\A^{ \text{s}}$.  

For every $(\varphi, m) \in  \A^{ \text{s}}$  satisfying 
\begin{equation}
\label{e:constraint-<}
\int_{\Om} \varphi \, \di x <V \qquad \text{and} \qquad \int_{Y} m \, \di y <W\,,
\end{equation}
 a recovery sequence can be constructed following~\cite{MR866718, MR0445362} for both~$\varphi$ and~$m$, relying on the convergence properties of~$S_{\eps}$ and~$S$. If one of the two constraints in~\eqref{e:constraint-<} is fulfilled with the equality, we may first construct a recovery sequence for the pair 
 \begin{displaymath}
 (\varphi_{\delta}, m_{\delta}) := ( \delta \varphi, \delta ( m-1) + 1) \quad \in \A^{ \text{s}}  \text{ for $\delta \in (0, 1)$}
 \end{displaymath}
 and proceed with a diagonal argument  as  $\delta \nearrow 1$.
\end{proof}

\subsection{Optimality conditions.} 
\label{s:opt-sharp}

We now focus on the  {\it Sharp-interface Microstructure-Topology
  Optimization} (SMTO)  problem 
\begin{align}
\label{e:min-sharp}
& \min_{\substack{\varphi \in BV(\Om; \{0, 1\})\,,\\ m \in  BV_{\#}
  (Y; \{1, 2\})}} \J^{ \text{s}}  (\varphi, m, u)  \ \ \ \text{subject to
  } \min_{u \in H^{1}_{\Gamma_{D}} (\Om; \R^{n})} \, \mathcal{E}
  (\varphi, m, u) \,. \tag{SMTO}
\end{align}
and in particular on its optimality conditions. From now on, we assume that the applied forces~$f$ and~$g$ belong to~$H^{1}(\Om; \R^{n})$ and to~$H^{1}(\Gamma_{N}; \R^{n})$, respectively. Such technical assumption is needed to state the optimality conditions in Theorem~\ref{t:optimality} (see also Proposition~\ref{p:differential}).

 We define the set of admissible variations w.r.t.~the phase-field variable~$\varphi$.

\begin{definition}[Admissible variations for $\varphi$]
\label{d:variations}
Let $\tau>0$. We say that $\Phi \in \mathrm{Adv}_{\tau}^{1}$ if the following holds:
\begin{itemize}
\item[$(\Phi_{0})$]  $\Phi \in C^{2,1}((-\tau, \tau) \times \overline{\Om}; \R^{n})$;

\item [$(\Phi_{1})$] there exists $C>0$ such that $\| \Phi(\cdot, x) - \Phi(\cdot, y) \|_{L^{\infty} (-\tau, \tau)} \leq C | x - y| $ for every $x, y \in \overline{\Om}$;

\item [$(\Phi_{2})$] $\Phi(t, x) \cdot \nu_{\Om} (x) = 0$ $\HH^{n-1}$-a.e.~on~$\partial \Om$, for every $t \in (-\tau, \tau)$, where $\nu_{\Om}$ denotes the outer unit normal vector to~$\partial\Om$;

\item [$(\Phi_{3})$] ${\rm spt}\Phi(t, \cdot) \cap \Gamma_{D} = \emptyset $ for every $t \in (-\tau, \tau)$.

\item[$(\Phi_{4})$] $\Phi(t, x) \cdot n_{\Gamma_{N}} (x)  = 0$ for every $x \in \partial_{\partial\Om} \Gamma_{N}$ and every $t \in (-\tau, \tau)$, where $\partial_{\partial\Om} \Gamma_{N}$ denotes the relative boundary of~$\Gamma_{N}$ in the topology of $\partial\Om$ and $n_{\Gamma_{N}}(x)$ stands for the outer unit normal vector to~$\Gamma_{N}$ in~$x$ laying in the tangent plane to~$\partial\Om$ in~$x$.   
\end{itemize}
The set of admissible variations for the phase-field variable 
$\varphi$  is $\mathrm{Adv}^{1} \coloneq \bigcup_{\tau >0} \mathrm{Adv}^{1}_{\tau}$.

 For $\Phi \in \mathrm{Adv}^{1}_{\tau}$, let $T_{t}$ be the
solution of the Cauchy problem 
\begin{equation}
\label{e:Tt}
\begin{cases}
\partial_{t} T_{t}(x) = \Phi(t, T_{t}(x)) &  \text{in $(-\tau, \tau)$,}
\\
T_{0}(x) = x \in \Om\,.
\end{cases}
\end{equation} 
 The map $T_t$  satisfies the following conditions:
\begin{itemize}
\item [$(T_{0})$] $T_{t}\colon \overline{\Om} \to \overline{\Om}$, $T_{t} (\Gamma_{N}) = \Gamma_{N}$, $T_{t} = {\rm Id}$ on~$\Gamma_{D}$, and $T \in C^{2,1}((-\tau, \tau) \times \overline{\Om} ; \R^{n})$;

\item[$(T_{1})$] $\| T_{t} - T_{0}\|_{C^{1,1}} \to 0$ as $t \to 0$;

\item [$(T_{2})$] $\partial_{t} T_{t}|_{t=0} = \Phi(0, \cdot)$, $\partial_{t} (\nabla_{x} T_{t})|_{t=0} = -  \partial_{t} (\nabla_{x} T_{t}) ^{-1}|_{t=0} = \nabla_{x} \Phi(0, \cdot)$, and $\partial_{t} (\det \nabla_{x} T_{t}) |_{t=0} = \dive \Phi(0, \cdot)$.
\end{itemize}
\end{definition}

Similarly, we introduce the variations w.r.t.~the microstructure variable~$m$.
\begin{definition}[Admissible variations for $m$]
\label{d:variations2}
Let $\tau>0$. We say that $\Psi \in \mathrm{Adv}_{\tau}^{2}$ if the following holds:
\begin{itemize}
\item[$(\Psi_{0})$]  $\Psi \in C^{2,1}((-\tau, \tau) \times Y; \R^{n})$ and $Y$-periodic;

\item [$(\Psi_{1})$] there exists $C>0$ such that $\| \Psi( \cdot, x ) - \Psi(\cdot, y) \|_{L^{\infty}(-\tau, \tau)} \leq C | x - y| $ for every $x, y \in Y$;

\item [$(\Psi_{2})$] $\Psi(t, y) \cdot \nu_{Y} (x) = 0$ $\HH^{n-1}$-a.e.~on~$\partial \Om$, for every $t \in (-\tau, \tau)$, where $\nu_{Y}$ denotes the external unit normal vector to~$\partial Y$.
\end{itemize}
The set of admissible variations for the microstructure variable 
$m$   is
$\mathrm{Adv}^{2} \coloneq \bigcup_{\tau >0} \mathrm{Adv}^2_{\tau}$.

 For $\Psi \in \mathrm{Adv}^{2}_{\tau}$ let $S_{t}$ be the
solution  of the Cauchy problem
\begin{equation}
\label{e:St}
\begin{cases}
\partial_{t} S_{t}(y) = \Psi(t, S_{t}(y)) &  \text{in $(-\tau, \tau)$,}
\\
S_{0}(y) = y \in Y\,.
\end{cases}
\end{equation} 
 The map $S_t$  satisfies the following conditions:
\begin{itemize}
\item [$(S_{0})$] $S_{t}\colon Y \to Y$ is $Y$-periodic, $S_{t} (\partial Y) = Y$, and $S \in C^{2,1}((-\tau, \tau) \times Y ; \R^{n})$;

\item[$(S_{1})$] $\| S_{t} - S_{0}\|_{C^{1,1}} \to 0$ as $t \to 0$;

\item [$(S_{2})$] $\partial_{t} S_{t}|_{t=0} = \Psi(0, \cdot)$, $\partial_{t} (\nabla_{y} S_{t})|_{t=0} = -  \partial_{t} (\nabla_{y} S_{t}) ^{-1}|_{t=0} = \nabla_{y} \Psi(0, \cdot)$, and $\partial_{t} (\det \nabla_{y} S_{t}) |_{t=0} = \dive \Psi (0, \cdot)$.
\end{itemize}
\end{definition}

%

Let us fix $(\varphi, m) \in BV(\Om; \{0, 1\}) \times BV_{\#} (Y; \{1, 2\})$ and $(\Phi, \Psi) \in {\rm Adv}^{1}_{\tau} \times {\rm Adv}_{\tau}^{2}$ for some $\tau>0$. For $t \in (-\tau, \tau)$, let us consider $T_{t} \colon \overline\Om \to \overline \Om $ and~$S_{t}\colon \overline{Y} \to \overline{Y}$ the maps defined through the Cauchy problems~\eqref{e:Tt} and~\eqref{e:St}, respectively. We set $(\varphi_{t}, m_{t}) \coloneq (\varphi \circ T_{t}^{-1}, m\circ S_{t}^{-1}) \in BV(\Om; \{0, 1\}) \times BV_{\#} (Y; \{1, 2\})$ and define, for $i, j=1, \ldots, n$ and $t \in (-\tau, \tau)$, $w_{i j}(t) (x, \cdot)  \in  H^{1}_{\#} (Y; \R^{n})$ as the unique solution to the correctors problem
\begin{align}
\label{e:T-corrector}
\min\bigg\{ \frac{1}{2} \int_{ Y} \C (\varphi_{t}, m_{t}) (x, y)  ( \strain_{y} (w) + e_{ij}) {\, \cdot\,} ( \strain_{y} (w) + e_{ij}) \, \di y  : & \, w \in H^{1}_{\#} (Y; \R^{n})
\\
&
\int_{Y} w(x, y) \, \di y = 0 \bigg\}. \nonumber
\end{align}
We further introduce the function $v_{ij} (t) \in L^{2}(\Om; H^{1}_{\#} (Y; \R^{n}))$ as
\begin{displaymath}
v_{i j} (t) (x, y) \coloneq w_{ij} (t) (T_{t}( x)  , S_{t} (y)) - \int_{Y} w_{ij} (t) (T_{t} (x)  , S_{t} (y))\, \di y \qquad \text{for $(x, y) \in \Om \times Y$.}
\end{displaymath}

The next two lemmas are devoted to the differentiability of $t \mapsto v_{i j} (t) $ in~$t = 0$.

\begin{lemma}[Bound on $v_{i j} (t)$]
\label{l:v-lip}
There exists a positive constant~$C=C(\Psi) > 0$ such that for every $i, j=1, \ldots, n$ and for a.e.~$x \in \Om$ the map $v_{ij} (x, \cdot) \colon (-\tau, \tau) \to  H^{1}_{\#} (Y; \R^{n})$ satisfies
\begin{equation}
\label{e:1001}
\| v_{ij} (t) - w_{ij}\|_{H^{1}( Y)} \leq C \big(  \| (\nabla_{y} S_{t}
)^{-1} - I \|_{L^{\infty} (Y) } + \| \det \nabla_{y} S_{t} - 1 \|_{L^{\infty} (Y) } \big)
\end{equation}
for every $t \in (-\tau, \tau)$.
\end{lemma}

\begin{proof}
We notice that, by Remark~\ref{r:0mean} and by the change of variables induced by~$S_{t}$, for a.e.~$x \in \Om$ the function~$v_{ij} (t) (x, \cdot)$ solves the minimum problem
\begin{align}
\label{e:mint}
& \min_{v \in H^{1}_{\#}(Y; \R^{n})}\bigg\{  \frac{1}{2}\! \int_{ Y} \! \C (\varphi_{t} , m ) (T_{t}
  (x) , y) ( \nabla_{y} v (\nabla_{ y} S_{t})^{-1} {+} e_{ij}) {\, \cdot\,}
  (\nabla_{y} v (\nabla_{ y} S_{t})^{-1} {+} e_{ij})\det \nabla_{ y} S_{t} \, \di y   \bigg\}
\\
& =  \min_{v \in H^{1}_{\#}(Y; \R^{n})}\bigg\{  \frac{1}{2}\! \int_{ Y} \! \C (\varphi , m ) (x , y) ( \nabla_{y} v (\nabla_{ y} S_{t})^{-1} + e_{ij}) {\, \cdot\,} (\nabla_{y} v (\nabla_{ y} S_{t})^{-1} + e_{ij})\det \nabla_{ y} S_{t} \, \di y\bigg\} \nonumber
\end{align}
From the minimum problem in~\eqref{e:mint}, from the fact that
$\int_{Y} v_{ij} (t) \, \di y = 0$, and from~$(S_1)$ of 
Definition~\ref{d:variations2}  we immediately deduce that
\begin{equation}
\label{e:1003}
\sup_{ \substack{i, j=1, \ldots, n, \\ t \in (-\tau, \tau)}} \|
v_{ij}( x,  t)\|_{H^{1}(Y)} < \infty
\end{equation} 
 for a.e. $x\in \Om$. 
Testing the Euler-Lagrange equations related to~\eqref{e:min-corrector} and~\eqref{e:mint} with~$v= v_{ij} (t) - w_{ij}$, subtracting one from the other, and rearranging the terms, we obtain that
\begin{align}
\label{e:1002}
0 &  =  \int_{Y} \C(\varphi, m) (\nabla_{y} v_{ij} (t) ( \nabla_{y} S_{t})^{-1} + e_{ij}) {\, \cdot\,} ( \nabla_{y} v_{ij} (t) - \nabla_{y} w_{ij}) ( \nabla_{y} S_{t})^{-1}  \, \det \nabla_{y} S_{t} \, \di y 
\\
&
\quad - \int_{Y} \C(\varphi, m) (\strain_{y} (w_{ij}) + e_{ij}) {\, \cdot\,} ( \strain_{y} (v_{ij} (t)) - \strain_{y} (w_{ij})) \, \di y \nonumber
\\
&
= \int_{Y} \C(\varphi, m) \big( \nabla_{y} v_{ij} (t) ( ( \nabla_{y} S_{t})^{-1} - I) \big) {\, \cdot\,} ( \nabla_{y} v_{ij} (t) - \nabla_{y} w_{ij}) ( \nabla_{y} S_{t})^{-1}  \, \det \nabla_{y} S_{t} \, \di y  \nonumber
\\
&
\quad + \int_{Y} \C(\varphi, m) (\strain_{y} (v_{ij} (t))  + e_{ij}) {\, \cdot\,} \big( ( \nabla_{y} v_{ij} (t) - \nabla_{y} w_{ij}) ( ( \nabla_{y} S_{t})^{-1} - I) \big)  \, \det \nabla_{y} S_{t} \, \di y \nonumber
\\
&
\quad + \int_{Y} \C(\varphi, m) (\strain_{y} (v_{ij} (t))+
  e_{ij}) {\, \cdot\,} ( \strain_{y} (v_{ij} (t)) -
  \strain_{y} ( w_{ij}))   \, ( \det \nabla_{y} S_{t} - 1)  \, \di y  \nonumber
\\
&
\quad + \int_{Y} \C(\varphi, m) (\strain_{y}( v_{ij} (t)) -
  \strain_{y} (w_{ij}) ) {\, \cdot\,} ( \strain_{y}(
  v_{ij} (t)) - \strain_{y} (w_{ij})) \, \di y \nonumber \,.
\end{align}
Hence, recalling~\eqref{e:hp-C} we infer from~\eqref{e:1002} that
\begin{align}
\label{e:1004}
\alpha & \| \strain_{y} v_{ij} (t) - \strain_{y}
         w_{ij}\|^{2}_{L^{2}(Y)} 
\\
&
\leq 2 \beta   \| \det \nabla_{y} S_{t}\|_{L^{\infty}(Y)} \| ( \nabla_{y} S_{t})^{-1}\|_{L^{\infty}(Y)} \big(\| v_{ij} (t) \|_{H^{1}(Y)}  + \| e_{ij}\|_{L^{\infty}(Y)} \big) \nonumber
\\
&
\qquad \cdot \| (\nabla_{y} S_{t})^{-1} - I\|_{L^{\infty}(Y)} \| \nabla_{y} v_{ij} (t) - \nabla_{y} w_{ij}\|_{L^{2}(Y)} \nonumber
\\
&
\qquad + \beta \big( \| v_{ij} (t) \|_{H^{1}(Y)}
  + \| e_{ij}\|_{L^{\infty}(Y)} \big)\,\|
  \strain_{y}( v_{ij} (t) )  - \strain_{y} (w_{ij})\|_{L^{2}(Y)} \| \det \nabla_{y} S_{t} - 1\|_{L^{\infty}(Y)}\,.  \nonumber
\end{align}
By Poincar\'e and Korn inequalities in $H^{1}_{\#} (Y; \R^{n})$ we deduce~\eqref{e:1001} from~\eqref{e:1003} and~\eqref{e:1004}.
\end{proof}

\begin{lemma}[Differentiability of $v_{ij}$]
\label{l:v-diff}
For every $i, j=1, \ldots, n$, the map $v_{ij} (x, \cdot) \colon
(-\tau, \tau) \to  H^{1}_{\#} (Y; \R^{n})$ is differentiable in~$t =
0$, uniformly w.r.t.~$x \in \Om$ (in the sense that the remainder can
be  uniformly estimated),  and its derivative $z_{ij}(x, \cdot) \coloneq \partial_{t}|_{t=0} v_{ij} (t)(x, \cdot) $ is the solution to
\begin{align}
\label{e:minz}
\min\bigg\{ \frac{1}{2}  \int_{Y} & \C(\varphi, m) (x, y)  \strain_{y} (z) {\, \cdot\,} \strain_{y} (z) \, \di y 
\\
&
\nonumber - \int_{Y} \C(\varphi, m) (x, y) \big( \nabla_{y} z\nabla_{y} \Psi(0) \big) {\,\cdot\,} ( \strain_{y}(  w_{ij} ) + e_{ij}) \, \di y 
\\
&
  - \int_{Y} \C(\varphi, m) (x, y)  \strain_{y} (z) {\,\cdot\,} \big( \nabla_{y} w_{ij} \nabla_{y} \Psi(0) \big) \, \di y \nonumber
  \\
  &
   + \int_{Y} \C(\varphi, m) (x, y) \strain_{y} (z) {\, \cdot\,}(  \strain_{y} (w_{ij}) + e_{ij}) \, \dive \Psi(0) \, \di y: \nonumber
   \\
   &
   \qquad \qquad \qquad z \in H^{1}_{\#} (Y; \R^{n}), \, \int_{Y} z \, \di y=0 \bigg\} \nonumber\,.
\end{align}
\end{lemma}

\begin{proof}
We proceed similarly to Proposition~\ref{p:C*-diff}, since we want to
show uniform differentiability in~$t=0$. For simplicity of notation,
we drop the dependence on~$x \in \Om$. We notice that a solution
to~\eqref{e:minz} exists and is unique, thanks to the  periodicity
 and the $0$-mean constraints.  We further notice that the $0$-mean condition in~\eqref{e:minz} can be dropped and $z_{ij}$ would still be a (non-unique) minimizer.

Writing the Euler-Lagrange equations related to~\eqref{e:min-corrector},~\eqref{e:mint}, and~\eqref{e:minz}, and rearranging the terms, for $t \in (-\tau, \tau)$ and $v \in H^1_{\#}(Y;\R^n)$ we have that
\begin{align}
\label{e:1000}
0  &=  \int_{Y} \C( \varphi, m)  ( \nabla_{y} v_{ij} (t)
     (\nabla_{ y} S_{t})^{-1} + e_{ij}) {\, \cdot\,} (\nabla_{y} v (\nabla_{ y} S_{t})^{-1} )\det \nabla_{ y} S_{t}\, \di y
\\
&\quad 
 - \int_{Y} \C(\varphi, m) ( \strain_{y} (w_{ij})  + e_{ij}) {\, \cdot\,} \nabla_{y} v\, \di y  - t \int_{Y} \C(\varphi, m) \strain_{y} (z_{ij}) {\, \cdot\,} \strain_{y} (v) \, \di y \nonumber
 \\
 &\quad
 + t \int_{Y} \C(\varphi, m) \big( \nabla_{y} v\nabla_{y} \Psi(0) \big) {\,\cdot\,} ( \strain_{y} ( w_{ij} ) + e_{ij}) \, \di y \nonumber
 \\
 &
 \quad
  + t \int_{Y} \C(\varphi, z)  \strain_{y} ( v ) {\,\cdot\,} \big( \nabla_{y} w_{ij} \nabla_{y} \Psi(0) \big) \, \di y \nonumber
  \\
  &\quad
   - t \int_{Y} \C(\varphi, m) \strain_{y} (v) {\, \cdot\,}(
    \strain_{y} (w_{ij}) + e_{ij}) \, \dive  \Psi(0) \, \di y\, \nonumber
   \\
   &
   = \int_{Y} \C( \varphi, m)  \big( \nabla_{y} v_{ij} (t)  ( (\nabla_{ y} S_{t})^{-1}- I + t\nabla_{ y}\Psi(0))  \big) {\, \cdot\,} (\nabla_{y} v (\nabla S_{t})^{-1} )\det \nabla_{ y} S_{t}\, \di y \nonumber
   \\
   &
   \quad +  \int_{Y} \C( \varphi, m)  ( \strain_{y}  (v_{ij} (t) )  + e_{ij}  ) {\, \cdot\,} \big(\nabla_{y} v ( (\nabla_{ y} S_{t})^{-1} ) - I + t\nabla_{ y}\Psi(0)) \big) \det \nabla_{ y} S_{t}\, \di y \nonumber
    \\
   &
   \quad +  \int_{Y} \C( \varphi, m)  ( \strain_{y}  ( v_{ij} (t) )  +
     e_{ij}  ) {\, \cdot\,} \strain_{y} (v) \big(  \det \nabla_{
     y} S_{t}- 1 - t\dive  \Psi(0) \big) \, \di y \nonumber
   \\
   &
   \quad + \int_{Y} \C(\varphi, m) \big( \strain_{y} (v_{ij} (t)) - \strain_{y} (w_{ij}) - t \strain_{y} (z_{ij}) \big) {\, \cdot\,} \strain_{y} (v) \, \di y \nonumber
   \\
   &
   \quad - t \int_{Y} \C(\varphi, m) \big( \nabla_{y} v_{ij}(t) \nabla_{y} \Psi(0) \big){\, \cdot\,} \big( \nabla_{y} v ( (\nabla_{y} S_{t})^{-1} - I) \big) \det \nabla_{y}S_{t} \, \di y   \nonumber
   \\
   &
    \quad - t \int_{Y} \C(\varphi, m) \big( \nabla_{y} v_{ij}(t) \nabla_{y}\Psi(0) \big){\, \cdot\,} \strain_{y} ( v )  \big(  \det \nabla_{y}S_{t} - 1 \big)\, \di y   \nonumber
    \\
    &
    \quad - t \int_{Y} \C(\varphi, m) \big( ( \nabla_{y} v_{ij} (t) - \nabla_{y} w_{ij}) \nabla_{y} \Psi(0) \big) {\, \cdot\,} \strain_{y} (v)\, \di y\nonumber
    \\
    &
    \quad - t \int_{Y} \C(\varphi, m) (\strain_{y} (v_{ij} (t))  + e_{ij}) {\, \cdot\,} (\nabla_{y} v \nabla_{y} \Psi (0)) (\det \nabla_{y} S_{t} - 1) \, \di y \nonumber
    \\
    &
    \quad - t \int_{Y} \C(\varphi, m) (\strain_{y}( v_{ij} (t)) - \strain_{y} ( w_{ij} ) ) {\, \cdot\,} (\nabla_{y} v \nabla_{y} \Psi (0)) \, \di y \nonumber
    \\
    &
    \quad + t \int_{Y} \C(\varphi, m) (\strain_{y} ( v_{ij} (t) ) -
      \strain_{y} ( w_{ij} ) ) {\, \cdot\,} \strain_{y} (v) 
      \, \dive  \Psi(0) \, \di y \nonumber\,.
\end{align}
Testing~\eqref{e:1000} with $v = v_{ij} (t) - w_{ij} - tz_{ij} \in H^{1}_{\#} (Y; \R^{n})$ and recalling~\eqref{e:hp-C}, we get that
\begin{align}
\label{e:2000}
\alpha  \| &\strain_{y}( v_{ij} (t) )   - \strain_{y} ( w_{ij} )  - t \strain_{y} (z_{ij})\|^{2}_{L^{2}(Y)}
\\
&
 \leq \Big( \beta \big(\| \nabla_{y} v_{ij} (t) \|_{L^{2}} + \| e_{ij}\|_{L^{\infty}(Y)} \big) \big(  \| (\nabla_{y} S_{t})^{-1} - I + t\nabla_{y} \Psi(0)\|_{L^{\infty}(Y)} \| \det \nabla_{y} S_{t}\|_{L^{\infty}(Y)} \big) \nonumber
 \\
 &
\qquad \quad \cdot \big( 1 + \| (\nabla S_{t})^{-1}\|_{L^{\infty} (Y)} \big)  \nonumber
 \\
 &
 \qquad +  \beta \big(\| \nabla_{y} v_{ij} (t) \|_{L^{2}(Y)} + \| e_{ij}\|_{L^{\infty}(Y)} \big) \| \det \nabla_{y} S_{t} - 1 - t\dive \Psi(0) \|_{L^{\infty}(Y)}  \big) \nonumber
  \\
  &
 \qquad  +  t\beta \| \nabla_{y} v_{ij} (t) \|_{L^{2}(Y)} \| \nabla_{y}
    \Psi(0)\|_{L^{\infty}(Y)}  \cdot \nonumber\\
&\qquad \quad \cdot \big(\| \det \nabla_{y} S_{t} -
    1\|_{L^{\infty}(Y)} + \| (\nabla_{y} S_{t})^{-1} -
    I\|_{L^{\infty}(Y)} \| \det \nabla S_{t}\|_{L^{\infty} (Y)}\big)
    \nonumber\\
&\qquad  + t\beta \big( \| \nabla_{y} v_{ij} (t) \|_{L^{2}(Y)} +
  \| e_{ij}\|_{L^{\infty} (Y)}\big)  \| \nabla_{y}
  \Psi(0)\|_{L^{\infty}(Y)}\cdot \nonumber \\
&\qquad \quad \cdot \| \det \nabla_{y} S_{t} -
    1\|_{L^{\infty}(Y)}\| \det \nabla_{y} S_{t} -
    1\|_{L^{\infty}(Y)}\nonumber
 \\ 
  &
 \qquad +   2  t \beta \big(  \| \nabla _{y} \Psi (0) \|_{L^{\infty}(Y)} \| \nabla_{y} v_{ij} (t) - \nabla_{y} w_{ij}\|_{L^{2}(Y)} \Big) \| \nabla_{y} v_{ij} (t)  - \nabla_{y}  w_{ij}  - t \nabla_{y} z_{ij}\|_{L^{2} (Y)} \,. \nonumber
\end{align}
Then, we infer from  Definition~\ref{d:variations2},  from Lemma~\ref{l:v-lip}, from~\eqref{e:2000}, and from Poincar\'e and Korn inequalities in~$H^{1}_{\#} (Y; \R^{n})$ that
\begin{displaymath}
\|  v_{ij} (t)  -  w_{ij}  - t  z_{ij}\|_{H^{1}(Y)} \leq {\rm o} (t)\,,
\end{displaymath}
where ${\rm o}  (t)$ does not depend on~$x \in \Om$ and is such that ${\rm o}  (t)/t\to 0$ as $t \to 0$. This concludes the proof of the lemma.
\end{proof}

\begin{corollary}[Differentiability of $\C^{*}(\varphi, m_{t} ) (x)$]
\label{c:C-diff-again}
The map $(-\tau, \tau) \ni t \mapsto \C^{*}(\varphi, m_{t} ) (x) \in \mathbb{R}^{n\times n\times n\times n}$ is differentiable in~$t=0$ uniformly w.r.t.~$x \in \Om$ and its derivative $\widetilde{\C}^{*} (\varphi, m) \coloneq \partial_{t}|_{t=0}  \C^{*}(\varphi, m_{t} ) $ is defined as
\begin{align}
\label{e:tildeC}
\widetilde{\C}^{*}_{ijkl} (\varphi, m) (x) = & \int_{Y} \C_{ijkl} (\varphi, m) (x, y) \big( \strain_{y}( z_{ij}) - \nabla_{y} w_{ij} \nabla_{y} \Psi(0) \big)  {\, \cdot\,} (\strain_{y}( w_{kl} )  + e_{kl}) \, \di y
\\
&
+ \int_{Y} \C_{ijkl} (\varphi, m) (x, y) ( \strain_{y} (w_{ij}) + e_{ij}) {\, \cdot\,} \big(  \strain_{y} (z_{kl}) - \nabla_{y} w_{kl} \nabla_{y} \Psi(0) \big) \, \di y \nonumber
\\
&
+ \int_{Y} \C_{ijkl} (\varphi, m) (x, y) ( \strain_{y} (w_{ij}) + e_{ij}) {\, \cdot\,} (\strain_{y} (w_{kl}) + e_{kl}) \dive \Psi(0) \, \di y\,,  \nonumber
\end{align}
where $z_{ij}= \partial_{t}|_{t=0} v_{ij} (t) \in H^{1}_{\#} (Y;
\R^{n})$ is the function  identified  in Lemma~\emph{\ref{l:v-diff}}.
\end{corollary}

\begin{proof}
For $t \in (-\tau, \tau)$ and~$x \in \Om$ we write $\C^{*} (\varphi, m_{t}) (x)$ using formula~\eqref{e:C-hom} and performing the change of variable induced by~$S_{t}$:
\begin{align}
\label{e:4000}
& \C^{*}_{ijkl}   (\varphi, m_{t}) (x)  = \int_{Y} \C(\varphi, m_{t}) (x, y) ( \strain_{y} (w_{ij} (t)) + e_{ij}) {\, \cdot\,} (\strain_{y} (w_{kl} (t) ) + e_{kl}) \, \di y
\\
&
=   \int_{Y} \C_{ijkl}(\varphi, m) (x, y) ( \nabla_{y} v_{ij} (t)  (\nabla_{y} S_{t})^{-1} + e_{ij}) {\, \cdot\,} (\nabla_{y} v_{kl} (t) (\nabla_{y} S_{t})^{-1}  + e_{kl}) \, \det \nabla_{y} S_{t} \, \di y .\nonumber
\end{align}
Then, the differentiability of~$\C^{*}_{ijkl}   (\varphi, m_{t}) (x) $
follows from the differentiability of~$t \mapsto v_{ij} (t)(x, \cdot)$
in~$t = 0$,  uniform w.r.t.~$x \in \Om$,  shown in
Lemma~\ref{l:v-diff}. Formula~\eqref{e:tildeC} also follows from
Lemma~\ref{l:v-diff}, together with the property of the map~$S_{t}$
reported in  Definition~\ref{d:variations2}.  
\end{proof}

Before stating the optimality conditions for~\eqref{e:min-sharp}, we prove the differentiability of the map~$t \mapsto u(t) \circ T_{t}$ in~$t=0$, where~$u(t) \in H^{1}_{\Gamma_{D}} (\Om; \R^{n})$ is the solution to
\begin{equation}\label{e:pb-f}
\min_{v \in H^{1}_{\Gamma_{D}} (\Om; \R^{n}) } \, \E(\varphi_{t} , m_{t}, v) \,.
\end{equation}

\begin{proposition}[Differentiability of $u(t)\circ T_t$]
\label{p:differential}
Let $f \in H^{1} (\Om; \R^{n})$, $g \in H^{1} (\Gamma_{N} ; \R^{n})$,
$(\Phi, \Psi) \in {\rm Adv}^{1} \times {\rm Adv}^{2}$, and~$T_{t},
S_{t}$ be as in 
Definitions~\emph{\ref{d:variations}-\ref{d:variations2}},  for $t \in (-\tau, \tau)$, respectively. Let $(\varphi, m) \in BV(\Om; \{0, 1\}) \times BV_{\#} (Y; \{1, 2\})$ and let $(\varphi_{t}, m_{t}) = (\varphi \circ T_{t}^{-1}, m\circ S_{t}^{-1})$. Then, the map $t \mapsto u(t) \circ T_{t}$ is differentiable in $t=0$ and $u_{(\Phi, \Psi)} \coloneq \partial_{t} ( u(t) \circ T_{t})|_{t=0} \in H^{1}_{\Gamma_{D}}(\Om; \R^{n})$ solves
\begin{align}
\label{e:opt-1}
\int_{\Om}  \C^{*}(\varphi, m) & \strain(u_{(\Phi, \Psi)} ) {\, \cdot\,} \strain (z) \, \di x 
+  \int_{\Om} \widetilde{\C}^{*}(\varphi, m) \strain( u ) {\, \cdot\,} \strain (z) \, \di x
\\
&
\nonumber - \int_{\Om} \C^{*} (\varphi, m) (\nabla u  \nabla \Phi(0) ){\, \cdot\,} \strain (z) \, \di x
 - \int_{\Om}\C^{*}(\varphi, m) \strain ( u )  {\, \cdot\,} (\nabla z \nabla \Phi(0) )\, \di x 
 \\
 &
 + \int_{\Om} \C^{*}(\varphi, m) \strain(u) {\, \cdot\,} \strain (z) \, \dive \Phi (0) \, \di x  - \int_{\Om} \varphi (\nabla f \Phi(0))  {\, \cdot\,}  z \, \di x \nonumber
 \\
 &
    - \int_{\Om} \varphi f {\, \cdot\,} z \, \dive \Phi(0)\, \di x - \int_{\Gamma_{N}} (\nabla g \Phi(0)) {\, \cdot\,} z \, \di \HH^{n-1} \nonumber
 \\
 &
 - \int_{\Gamma_{N}} g {\, \cdot\,} z \, \big(  \dive \Phi(0) - \nu_{\Om} (x) \cdot \nabla \Phi(0) \nu_{\Om} (x) \big)  \, \di \HH^{n-1}=0\nonumber
\end{align}
for every $z \in H^{1}_{\Gamma_{D}} (\Om; \R^{n})$.
\end{proposition}

In the proof of Proposition~\ref{p:differential} we will use the following lemma, whose proof can be found in, e.g.,~\cite[Lemma~4.1]{KneMie2008} or~\cite[Lemma~3.8]{almi2015}.

\begin{lemma}
\label{l:1}
For every $v \in H^{1}(\Om; \R^{n})$ and every $\Phi \in
\mathrm{Adv}^{1}$ we have that the quotient $t^{-1} (v - v \circ
T_{t}^{-1})$ converges to $\nabla v \Phi(0)$  strongly  in $L^{2}(\Om; \R^{n})$ as $t \to 0$.
\end{lemma}

\begin{proof}[Proof of Proposition~\ref{p:differential}]
Let $\tau>0$ be such that $(\Phi, \Psi) \in  \mathrm{Adv}^{1}_{\tau} \times {\rm Adv}^{2}_{\tau}$. 
By the change of coordinates $x = T_{t} (z) $ we get that the function $u(t) \circ T_{t}$ is solution to
\begin{align}
\label{e:3}
\min  \bigg\{ \frac{1}{2} & \int_{\Om} \C^{*}(\varphi, m_{t}) (\nabla v (\nabla T_{t})^{-1})  {\, \cdot\,} (\nabla v (\nabla T_{t})^{-1} )\, \det \nabla T_{t} \, \di x  - \int_{\Om} \varphi (f \circ T_{t})  {\, \cdot\,} v \, \det \nabla T_{t}   \, \di x
\\
&
 - \int_{\Gamma_{N}} (g \circ T_{t}) {\, \cdot\,} v \,\det \nabla T_{t} |  \nabla T_{t} (x) \nu_{\Om} ( x )   |  \, \di \HH^{d-1}: \, v \in H^{1}_{\Gamma_{D}}( \Om; \R^{n})\bigg\} \,. \nonumber
\end{align}
We define the map $F(\cdot, \cdot) \colon (-\tau, \tau) \times H^{1}_{\Gamma_{D}}(\Om; \R^{n}) \to (H^{1}_{\Gamma_{D}} (\Om; \R^{n}) )'$ as 
\begin{align}
\label{e:4}
F(t, \psi)[\zeta] \coloneq &\int_{\Om} \C^{*} (\varphi, m_{t}) (\nabla \psi (\nabla T_{t})^{-1}) {\, \cdot\,} ( \nabla \zeta ( \nabla T_{t})^{-1}) \det \nabla T_{t} \, \di x 
 \\
 &
 \nonumber - \int_{\Om} \varphi (f \circ T_{t})  {\, \cdot\,} \zeta  \det \nabla T_{t} \, \di x 
 - \int_{\Gamma_{N}} (g \circ T_{t}) {\, \cdot\,} \zeta \,\det \nabla T_{t} |  \nabla T_{t} (x) \nu_{\Om}  (x)   | \, \di \HH^{d-1}  \nonumber
\end{align}
for $t \in (-\tau, \tau)$ and $\psi, \zeta \in H_{\Gamma_{D}}^{1}(\Om;
\R^{d})$. Then, $F$ satisfies $F(t, u(t) \circ T_{t}) = 0$, as $u(t)
\circ T_{t}$ solves~\eqref{e:3} and the map $t \mapsto u(t) \circ
T_{t}$ is continuous. Thanks to Corollary~\ref{c:C-diff-again} and to
the regularity of~$f$ and of~$g$, the map $t \mapsto F(t, \psi)$ is
differentiable in $t=0$ uniformly w.r.t.~$\psi \in H^{1}_{\Gamma_{D}}
(\Om; \R^{n})$ such that $\|\psi\|_{H^{1}} \leq R$, while $\psi
\mapsto F(t, \psi)$ is Lipschitz continuous uniformly in $t \in
(-\tau, \tau)$ by 
Definitions~\ref{d:variations}-\ref{d:variations2}.  Moreover, we have that for every $\eta, \zeta \in H^{1}_{\Gamma_{D}}( \Om; \R^{n}) $ and every $t \in (-\tau, \tau)$
\begin{align*}
\big(D_{\psi} F(t, v) [\zeta]\big) [\eta] = \int_{\Om} \C^{*}(\varphi, m_{t}) (\nabla \eta (\nabla T_{t})^{-1}) {\, \cdot\,} (\nabla \zeta (\nabla T_{t})^{-1}) \det \nabla T_{t} \, \di x\,.
\end{align*} 
Thus, we deduce from Lemma~\ref{l:v-lip}, from
equality~\eqref{e:4000}, and from 
Definitions~\ref{d:variations}-\ref{d:variations2}  that
\begin{align*}
\| D_{\psi}  F(t, u) [\zeta] & - D_{\psi} F(0, u) [\zeta]\|_{(H^{1}_{\Gamma_{D}} (\Om; \R^{n}) )'}
\\
&
 \leq C (\| \mathrm{I} - (\nabla T_{t}) ^{-1} \|_{ L^\infty(\Om)} + \| 1 - \det \nabla T_{t}\|_{ L^\infty(\Om)} ) \| \zeta\|_{ H^{1}_{\Gamma_{D}} (\Om)}
 \\
 &
 \qquad+  ( \| (\nabla_{y} S_{t} )^{-1} - I \|_{L^{\infty} (Y) } + \| \det \nabla_{y} S_{t} - 1 \|_{L^{\infty} (Y) }  ) \| \zeta\|_{ H^{1}_{\Gamma_{D}} (\Om)}\,.
\end{align*}
and that $D_{\psi} F(0, v) \colon H^{1}_{\Gamma_{D}}( \Om; \R^{n}) \to
(H^{1}_{\Gamma_{D}, 0} ( \Om ; \R^{n}))'$ is invertible
(isomorphism). Hence, we can apply the Implicit Function Theorem 
in  \cite[Theorem~3.1]{Ulbrich} (see also~\cite{Simon, Wach}) and
 obtain  the differentiability of the map $t \mapsto u (t) \circ T_{t}$ in $t=0$. Denoting $u_{(\Phi, \Psi)} \coloneq \partial_{t}( u(t) \circ T_{t}) |_{t=0} \in H^{1}_{\Gamma_{D}} (\Om; \R^{n})$, we further have
\begin{displaymath}
D_{\psi} F(0, u) [u_{(\Phi, \Psi)}] = - \partial_{t} F (0, u)\,.
\end{displaymath}
By Lemma~\ref{l:v-diff}, by Corollary~\ref{c:C-diff-again},
by~\cite[Theorem~7.31]{APF}, and by 
Definitions~\ref{d:variations}-\ref{d:variations2},  the last equality may be rewritten as
\begin{align*}
 &\int_{\Om} \C^{*}( \varphi, m) \strain (u_{(\Phi, \Psi)}) {\,
   \cdot\,} \strain (z) \, \di x \\
&\quad = - \int_{\Om} \widetilde{\C}^{*}(\varphi, m) \strain( u ) {\, \cdot\,} \strain (z) \, \di x 
 + \int_{\Om} \C^{*}(\varphi, m) (\nabla u \nabla \Phi(0))  {\,
  \cdot\,} \strain ( z)  \, \di x 
\\
 &\qquad
 + \int_{\Om} \C^{*}(\varphi, m) \strain (u) {\, \cdot\,} (\nabla z \nabla \Phi(0)) \, \di x 
 - \int_{\Om} \C^{*}(\varphi, m) \strain (u) {\, \cdot\,} \strain (z)\, \dive \Phi(0) \, \di x 
 \\
 &\qquad+ \int_{\Om}\varphi (\nabla  f \Phi(0)){\, \cdot\,}  z  \, \di x + \int_{\Om} \varphi f {\, \cdot\,} z \, \dive \Phi(0)\, \di x  
 \\
 &\qquad
 + \int_{\Gamma_{N}} (\nabla g \Phi(0)) {\, \cdot\,} z \, \di
   \HH^{n-1}  
 + \int_{\Gamma_{N}} g {\, \cdot\,} z \, \big(  \dive \Phi(0) - \nu_{\Om} (x) \cdot \nabla \Phi(0) \nu_{\Om} (x) \big)  \, \di \HH^{n-1} \,,
\end{align*}
for every $z \in H^{1}_{\Gamma_{D}}( \Om; \R^{n})$. This concludes the proof of the proposition.
\end{proof}

\begin{corollary}[Differentiability of $\J^{ \text{s}}(\varphi_{t}, m_{t}, u(t))$]
\label{c:I-diff}
Under the assumptions of Proposition~\emph{\ref{p:differential}}, the map $t \mapsto \J^{ \text{s}}  (\varphi_{t}, m_{t}, u(t))$ is differentiable in~$t=0$ and its derivative writes
\begin{align}
\label{e:I-derivative}
\partial_{t} &\J^{ \text{s}}  (\varphi_{t}  , m_{t}, u(t))|_{t=0}  =   - \int_{\Om} \widetilde{\C}^{*} (\varphi, m) \strain( u ){\, \cdot\,} \strain (u) \, \di x
\\
& 
 +2 \int_{\Om}\C^{*}(\varphi, m) \strain ( u )  {\, \cdot\,} (\nabla u \nabla \Phi(0))\, \di x  
 -  \int_{\Om} \C^{*}(\varphi, m) \strain(u) {\, \cdot\,} \strain (u) \, \dive \Phi (0) \, \di x  \nonumber
\\
& 
 + 2\int_{\Om}\varphi (\nabla  f \Phi(0)){\, \cdot\,}  u  \, \di x + 2 \int_{\Om} \varphi f {\, \cdot\,} u \, \dive \Phi(0)\, \di x  
 + 2\int_{\Gamma_{N}} (\nabla g \Phi(0)) {\, \cdot\,} u \, \di \HH^{n-1} \nonumber
 \\
 &
  + 2 \int_{\Gamma_{N}} g {\, \cdot\,} u \, \big(  \dive \Phi(0) - \nu_{\Om} (x) \cdot \nabla \Phi(0) \nu_{\Om} (x) \big)  \, \di \HH^{n-1} \nonumber
  \\
  &
+ c_{H}\!\! \int_{\Om} ( \dive \Phi(0) - \nu_{\varphi} \cdot \nabla \Phi(0) \nu_{\varphi}) \, \di | D \varphi| 
 + c_{H}\!\! \int_{Y} ( \dive \Psi(0) - \nu_{m} \cdot \nabla \Psi(0) \nu_{m}) \, \di | D m|\,,  \nonumber
\end{align}
where $\nu_{\varphi}$ and~$\nu_{m}$ denote the approximate interior unit normal to~$\{\varphi  = 1\}$ and to~$\{m = 2\}$, respectively, and~$\widetilde{\C}^{*}$ is the tensor defined in~\eqref{e:tildeC}.
\end{corollary}

\begin{proof}
Let us fix $(\Phi, \Psi) \in \mathrm{Adv}^{1} \times {\rm Adv}^{2}$. By Proposition~\ref{p:differential} the function $t \mapsto u(t) \circ T_{t}$ with~$u(t)$ solution to~\eqref{e:pb-f} is differentiable in~$t=0$ with derivative~$u_{(\Phi, \Psi) } \coloneq \partial_{t} (u(t) \circ T_{t})|_{t=0} \in H^{1}_{\Gamma_{D}} (\Om; \R^{n})$ satisfying~\eqref{e:opt-1}. Combining~\eqref{e:opt-1} and the minimality of~$u$ we obtain
\begin{align}
\label{e:1}
\int_{\Om} \varphi f {\, \cdot\,} & u_{(\Phi, \Psi)} \, \di x   + \int_{\Gamma_{N}} g{\, \cdot\,} u_{(\Phi, \Psi)}\, \di \HH^{n-1}  =  \int_{\Om} \C^{*}(\varphi, m) \strain (u) {\, \cdot\,} \strain (u_{(\Phi, \Psi)}) \, \di x
\\
&
 =  - \int_{\Om} \widetilde{\C}^{*}(\varphi, m) \strain(u ) {\, \cdot\,} \strain( u) \, \di x 
 + 2 \int_{\Om} \C^{*}(\varphi, m) \strain (u) {\, \cdot\,} (\nabla u \nabla \Phi(0)) \, \di x \nonumber
 \\
 &
 \qquad - \int_{\Om} \C^{*}(\varphi, m) \strain (u) {\, \cdot\,} \strain (u) \dive \Phi(0) \, \di x 
 + \int_{\Om}\varphi (\nabla  f \Phi(0)){\, \cdot\,}  u  \, \di x  \nonumber
 \\
 &
\qquad + \int_{\Om} \varphi f {\, \cdot\,} u \, \dive \Phi(0)\, \di x  
 + \int_{\Gamma_{N}} (\nabla g \Phi(0)) {\, \cdot\,} u \, \di \HH^{n-1} \nonumber
 \\
 &
 \qquad + \int_{\Gamma_{N}} g {\, \cdot\,} u \, \big(  \dive \Phi(0) - \nu_{\Om} (x) \cdot \nabla \Phi(0) \nu_{\Om} (x) \big)  \, \di \HH^{n-1} \,.\nonumber
 \end{align}

We now compute the derivative of cost functional~$\J^{ \text{s}} $. For
$|t|$ small enough, using the notations of
Proposition~\ref{p:differential},  we have that 
\begin{align}
\label{e:2}
 &\J^{ \text{s}}  (\varphi_{t}, m_{t}, u(t))  - \J^{ \text{s}}  (\varphi, m, u) 
\\
&
\quad= \int_{\Om} \varphi (f \circ T_{t} ) {\, \cdot\,} (u(t) \circ
  T_{t}) \det \nabla T_{t}\,\di x \nonumber\\
&\qquad + \int_{\Gamma_{N}} \!\!\!\! (g\circ T_{t}) {\, \cdot \,} (u(t) \circ T_{t}) \det \nabla T_{t} | (\nabla T_{t})^{-1} \nu_{\Om} | \, \di \HH^{n-1} \nonumber
\\
&
\qquad + c_{ H} \big( \mathcal{P} (\{\varphi_{t}=1\}; \Om) + \mathcal{P} (\{ m_{t}=2\}; Y) \big)  \nonumber
\\
&
\qquad -  \int_{\Om} \varphi f {\, \cdot\,} u \, \di x
  -\int_{\Gamma_{N}} g {\, \cdot\,} u \, \di \HH^{n-1}  -  c_{ H} \big( \mathcal{P} (\{\varphi =1\}; \Om) + \mathcal{P} (\{ m =2\}; Y) \big)    \nonumber
\\
&
\quad = \int_{\Om} \varphi ( f \circ T_{t}) {\, \cdot\,} ( u(t) \circ T_{t} - u ) \, \det \nabla T_{t} \, \di x + \int_{\Om} \varphi ( f \circ T_{t})    {\, \cdot\,} u (\det \nabla T_{t} - 1)  \, \di x   \nonumber
\\
&
 \qquad + \int_{\Om} \varphi ( f \circ T_{t} - f) {\,\cdot\,} u \, \di x  + \int_{\Gamma_{N}} ( (g\circ T_{t}) - g) \cdot (u(t) \circ T_{t}) \det \nabla T_{t} | (\nabla T_{t})^{-1} \nu_{\Om} | \, \di \HH^{n-1} \nonumber
 \\
 &
\qquad + \int_{\Gamma_{N}} g {\, \cdot\,} (u(t) \circ T_{t} - u) \det
   \nabla T_{t} | (\nabla T_{t})^{-1} \nu_{\Om} | \, \di \HH^{n-1}
   \nonumber\\
&\qquad +  \int_{\Gamma_{N}} g {\, \cdot\,}  u \big(  \det \nabla T_{t} | (\nabla T_{t})^{-1} \nu_{\Om} |  - 1\big) \, \di \HH^{n-1} \nonumber
 \\
 &
\qquad + c_{H} \big(  \mathcal{P} (\{\varphi_{t}=1\}; \Om) - \mathcal{P} (\{\varphi =1\}; \Om) \big)  \vphantom{\int_{\Om} }+ c_{H} \big(   \mathcal{P} (\{ m_{t}=2\}; Y) -  \mathcal{P} (\{ m =2\}; Y) \big)\,.  \nonumber
\end{align}
We now divide~\eqref{e:2} by $t$ and pass to the limit as~$t \to
0$. In view of  Definition~\ref{d:variations}  and of Proposition~\ref{p:differential} we obtain
\begin{align}
\label{e:3.1}
\partial_{t}  \J^{ \text{s}}  (\varphi_{t}  , m_{t}, u(t))|_{t=0}  = &  \int_{\Om} \varphi f {\, \cdot\,} u_{(\Phi, \Psi)} \, \di x  + \int_{\Om} \varphi f{\, \cdot\,} u \, \dive \Phi(0) \, \di x + \int_{\Om} \varphi ( \nabla f \Phi(0)) {\, \cdot\,} u \, \di x 
\\
& 
+ \int_{\Gamma_{N}} (\nabla g \Phi(0)) {\, \cdot\,} u \, \di \HH^{n-1} + \int_{\Gamma_{N}} g {\, \cdot\,} u_{(\Phi, \Psi)} \, \di \HH^{n-1} \nonumber
\\
&
+ \int_{\Gamma_{N}} g {\, \cdot\,} u \big( \dive \Phi(0) - \nu_{\Om} {\, \cdot\,} \nabla \Phi(0) \nu_{\Om} \big) \, \di \HH^{n-1} \nonumber
\\
&
+ c_{H} \int_{\Om} ( \dive \Phi(0) - \nu_{\varphi} \cdot \nabla \Phi(0) \nu_{\varphi}) \, \di | D \varphi|  \nonumber
\\
&
 + c_{H} \int_{Y} ( \dive \Psi(0) - \nu_{m} \cdot \nabla \Psi(0) \nu_{m}) \, \di | D m|\,,  \nonumber
\end{align}
where we have used the formula for the first variation of the perimeter (see, e.g., \cite[Section~7.3]{APF}). Inserting~\eqref{e:1} into~\eqref{e:3.1} we infer~\eqref{e:opt-3}.
\end{proof}

We are now in a position to state the optimality conditions for~\eqref{e:min-sharp}.

\begin{theorem}[Optimality conditions for \eqref{e:min-sharp}]
\label{t:optimality}
Let 
$f \in
H^{1}(\Om; \R^{n})$, $g \in H^{1} (\Gamma_{N}; \R^{n})$, and
$(\varphi, m) \in \mathcal{A}^{ \text{s}} $ be a solution
to~\eqref{e:min-sharp} with corresponding displacement $u \in
H^{1}_{\Gamma_{D}} (\Om; \R^{n})$. Then, there exist~$(\lambda, \mu)
\in \R^{2}$ such that for every $(\Phi, \Psi) \in  {\rm Adv}^{1}
\times {\rm Adv}^{2}$  the following holds:
\begin{align}
\label{e:opt-3}
0 = & - \int_{\Om} \widetilde{\C}^{*} (\varphi, m) \strain( u  ){\, \cdot\,} \strain (u) \, \di x + 2\int_{\Om} \C^{*} (\varphi, m) (\nabla u  \nabla \Phi(0) ){\, \cdot\,} \strain (u) \, \di x
  \\
 &
 - \int_{\Om} \C^{*}(\varphi, m) \strain(u) {\, \cdot\,} \strain (u) \, \dive \Phi (0) \, \di x  + 2\int_{\Om}\varphi (\nabla  f \Phi(0)){\, \cdot\,}  u  \, \di x \nonumber
 \\
 &
 + 2 \int_{\Om} \varphi f {\, \cdot\,} u \, \dive \Phi(0)\, \di x  
 + 2\int_{\Gamma_{N}} (\nabla g \Phi(0)) {\, \cdot\,} u \, \di \HH^{n-1} \nonumber
 \\
 &
  + 2 \int_{\Gamma_{N}} g {\, \cdot\,} u \, \big(  \dive \Phi(0) - \nu_{\Om} (x) \cdot \nabla \Phi(0) \nu_{\Om} (x) \big)  \, \di \HH^{n-1} \nonumber \nonumber
\\
& 
+ c_{W}\!\! \int_{\Om} ( \dive \Phi(0) - \nu_{\varphi} \cdot \nabla \Phi(0) \nu_{\varphi}) \, \di | D \varphi| 
 + c_{W}\!\! \int_{Y} ( \dive \Psi(0) - \nu_{m} \cdot \nabla \Psi(0) \nu_{m}) \, \di | D m| \nonumber
 \\
 &
 + \lambda \int_{\Om} \varphi \, \dive \Phi(0) \, \di x + \mu \int_{Y} m \, \dive \Psi(0)\, \di y \,, \nonumber
\end{align}
where~$\widetilde{\C}^{*}$ is the tensor defined in~\eqref{e:tildeC}.
\end{theorem}

\begin{proof}
We follow the lines of~\cite[Theorem~13]{Garcke_16}. Let $(\varphi, m)
\in  \mathcal{A}^{ \text{s}} $ be a solution to~\eqref{e:min-sharp},
let us fix $(\Phi, \Psi) \in \mathrm{Adv}^{1} \times {\rm
  Adv}^{2}$,  and let us denote by~$(T_{t}, S_{t})$ the fields defined
in  Definitions~\ref{d:variations}-\ref{d:variations2}.  If 
\begin{displaymath}
\int_{\Om} \varphi \, \di x < V \qquad \text{and} \qquad \int_{Y} m \, \di y <W\,,
\end{displaymath}
then equality~\eqref{e:opt-3} follows from Corollary~\ref{c:I-diff} with $(\lambda, \mu) = (0,0)$, since $(\varphi_{t}, m_{t}) \coloneq (\varphi\circ T_{t}^{-1}, m \circ S_{t}^{-1}) \in \mathcal{A}^{ \text{s}} $ for $| t| \ll 1$. 

Let us assume that $\int_{\Om}\varphi \, \di x = V$ and $\int_{Y} m \,
\di y =W$ (the other cases can be treated similarly). We fix $\Xi
\in {\rm Adv}^{1}$  and $\Theta \in  {\rm Adv}^{2}$  such that
\begin{displaymath}
\int_{\Om} \varphi \, \dive \Xi(0)\, \di x = 1 \qquad \text{and} \qquad \int_{Y} m \, \dive \Theta(0) \, \di y = 1\,,
\end{displaymath}
and consider the corresponding maps $M_{s}$ and~$N_{\sigma}$ solutions of
\begin{displaymath}
\left\{\begin{array}{ll}
\partial_{s} M_{s} (x) = \Xi(s, M_{s}(x))\,,\\[1mm]
M_{0} (x) = x \in \overline{\Om}\,,
\end{array} \right.
\qquad \text{and}\qquad
\left\{\begin{array}{ll}
\partial_{\sigma} N_{\sigma} (y) = \Theta(\sigma, N_{\sigma}(y))\,,\\[1mm]
N_{0} (y) = y \in Y\,,
\end{array} \right.
\end{displaymath}
defined for $|s|, \, |\sigma| \ll 1$. The map 
\begin{displaymath}
\eta (t, s, \sigma) \coloneq \bigg( \int_{\Om} \varphi \circ (M_{s} \circ T_{t})^{-1}\, \di x - V , \int_{Y} m \circ (N_{\sigma} \circ S_{t})^{-1} \, \di y - W \bigg)
\end{displaymath}
satisfies $\eta(0, 0, 0) = (0, 0)$ and 
\begin{displaymath}
\big( \partial_{s} \eta (t, s, \sigma) , \partial_{\sigma} \eta(t, s, \sigma) \big)|_{(t, s, \sigma) = (0,0,0)} =  \bigg( \int_{\Om} \varphi \, \dive \Xi(0)\, \di x \, ,  \int_{Y} m \, \dive \Theta(0) \, \di y\bigg) = (1, 1)\,.
\end{displaymath}
Then, by the Implicit Function Theorem there exist $t_{0}>0$ and a
curve $(s, \sigma) \colon (-t_{0}, t_{0}) \to \R^{2}$ such that~$(s,
\sigma)(-t_{0}, t_{0})$ is  in  a neighborhood of~$(0,0)$ and $\eta(t, s(t), \sigma(t)) = (0,0,0)$ for every $t \in (-t_{0}, t_{0})$. Furthermore, it holds 
\begin{equation}
\label{e:ssigma}
(s'(0), \sigma'(0)) = \partial_{t} \eta(t, s, \sigma) |_{(t, s, \sigma) = (0,0,0)} = \bigg( \int_{\Om} \varphi \, \dive \Phi(0)\, \di x \, ,  \int_{Y} m \, \dive \Psi(0) \, \di y\bigg)  \,.
\end{equation}
Then, the optimality of~$(\varphi, m)$ implies that
\begin{align}
\label{e:9000}
0 & = \partial_{t} \J^{ \text{s}} \big(\varphi \circ (M_{s(t)} \circ T_{t})^{-1} \, , m \circ (N_{\sigma(t)} \circ S_{t})^{-1} \big)|_{t=0}
\\
&
= \partial_{s}\J^{ \text{s}}  (\varphi \circ M_{s}^{-1}, m)|_{s=0} s'(0) + \partial_{\sigma} \J^{ \text{s}}  ( \varphi, m \circ N_{\sigma}^{-1})|_{\sigma=0} \sigma'(0) + \partial_{t}\J^{ \text{s}}  (\varphi \circ T_{t}^{-1}, m\circ S_{t}^{-1})|_{t=0}\,.  \nonumber
\end{align}
By Proposition~\ref{p:differential}, the first two derivatives on the right-hand side of~\eqref{e:9000} exist. Setting
\begin{displaymath}
(\lambda, \mu) \coloneq \Big( \partial_{s} \J^{ \text{s}}  (\varphi \circ M_{s}^{-1}, m)|_{s=0} \, ,  \partial_{\sigma} \J^{ \text{s}}  ( \varphi, m \circ N_{\sigma}^{-1})|_{\sigma=0} \Big)\,,
\end{displaymath}
we deduce~\eqref{e:opt-3} from~\eqref{e:ssigma}-\eqref{e:9000}. 
\end{proof}

\bigskip
\noindent\textbf{Acknowledgments}
The work of SA was partially funded by the Austrian Science Fund
through the projects ESP-61 and P-35359. SA also acknowledges the warm
hospitality of the University of Vienna, of the TU Wien, and of ESI
during the workshop \emph{Between Regularity and Defects: Variational
  and Geometrical Methods in Materials Science}, where part of this
research was carried out.  US is partially funded by the  Austrian
Science Fund grants I-5149, F-65, I-4354, \linebreak and~{P-32788}.

\bibliographystyle{siam}

\end{document}